\documentclass[preprint,3p,12pt]{elsarticle}
\usepackage{epsfig,latexsym,amssymb, amsmath,amscd,amsthm,amsxtra}
\usepackage{color}
\usepackage{graphicx}
\usepackage{CJK}
\usepackage{indentfirst}
\usepackage{setspace}
\usepackage{amsmath}
\usepackage{amssymb}
\usepackage{graphicx}
\usepackage{subfig}
\usepackage{enumerate}
\usepackage{float}
\usepackage{animate}
\usepackage{booktabs}
\usepackage{threeparttable}
\usepackage{multirow}
\theoremstyle{remark}
\newtheorem{remark}{\textbf{Remark}}[section]

 \numberwithin{equation}{section}

\journal{Elsevier}

\begin{document}
	
\begin{frontmatter}	

\title{The Fictitious Domain Method Based on Navier Slip Boundary Condition for Simulation of Flow-Particle Interaction }

\author[CAEP]{Rong Zhang}
\author[SCU]{Qiaolin He}
\address[CAEP]{Key Lab of Intelligent Analysis and Decision on Complex Systems,
Chongqing University of Posts and Telecommunications, Chongqing 400065, China, (zhangrong@cqupt.edu.cn)
}
\address[SCU]{School of Mathematics, Sichuan University, Chengdu, 610064, China, (qlhejenny@scu.edu.cn)}

\begin{abstract}
	In this article, we develop a least--squares/fictitious domain method for direct simulation of fluid particle motion  with Navier slip boundary condition at the fluid--particle interface.  Let $\Omega$ and $B$ be two bounded domains of $\mathbb{R}^{d}$ such that  $\overline{B} \subset \Omega$.  The motion of solid particle $B$ is governed by Newton's equations.  Our goal here is to develop a fictitious domain method where one solves a variant of the original problem on the full $\Omega$, followed by a well--chosen correction over $B$ and corrections related to translation velocity and angular velocity of the particle.  This method is of the virtual control type and relies on a least--squares formulation making the problem solvable by a conjugate gradient algorithm operating in a well chosen control space.  Since the fully explicit scheme to update the particle motion using Newton's equation is unstable, we propose and implement an explicit--implicit scheme in which, at each time step,  the position of the particle is updated explicitly, and the solution of Navier-Stokes equations and particle velocities are solved by the the  least--squares/fictitious domain method implicitly.  Numerical results are given to verify our numerical method.
\end{abstract}


\end{frontmatter}



\section{Introduction}
In order to understand the interactions between the solid particle and the fluid, we take a more fundamental approach with the Navier--Stokes equations  for fluid and Newton's equations of motion for the particle. Numerical simulation  of these equations can give us all the details of the flow and the fluid--particle interactions together with clear understandings of the mechanisms involved.   There are some numerical works on these types of  simulations.  Authors in \cite{Hu1} and \cite{Hu2} used arbitrary Lagrangian--Eulerian method to directly simulate the motion of fluid and particles, which  based on moving unstructured grids and needs remeshing and projection.

 Another  important method is fictitious domain method. In \cite{Glowinski1,Glowinski3,Glowinski2,Glowinski4}, fictitious domain methods were discussed  for the solution of Dirichlet problems, where the Dirichlet boundary condition being  enforced as a side constraint, using a {\em boundary supported Lagrange multiplier}. A volume--supported  Lagrange multiplier based fictitious domain method was introduced  in \cite{Glowinski5}, the main motivation being the direct numerical simulation of particulate flow when the number of particles  exceeds $10^{3}$. Initially tested on particulate flow with spherical particles, the method discussed in  \cite{Glowinski5} was generalized to situations involving particles with more complicated shapes, as shown for example in \cite{Glowinski6}. A (brief) history of {\em fictitious domain methods } can be found in, e.g., Chapter 8 in \cite{Handbook}.  Examples of non Lagrange multiplier based fictitious domain methods can be found in the {\em immersed boundary method} for the simulation of incompressible viscous flow in regions with elastic moving boundaries \cite{Peskin}. Then, M. Uhlmann \cite{Uhlmann} used the immersed boundary method to simulate particulate flow.  

The main idea behind fictitious domain methods is to extend a problem initially posed on a geometrically complex shaped domain to a larger simpler domain; this provides two main advantages  when constructing numerical schemes: (i) The extended domain is geometrically simpler and allows the use of fast solvers. (ii) The same fixed mesh can be used for the entire computation, eliminating thus the need for repeated remeshing and projection. All the studies that we know, concerning the application of fictitious domain methods to the simulation of particulate flow, consider no-slip boundary  conditions at the interface between fluid and particles.  There are situations however, in micro-fluidics for example, where a slip  condition on the particle surface is more realistic than the no-slip one.  If  the no-slip boundary condition on the particle surface is replaced by the Navier slip boundary condition, the volume--supported  Lagrange multiplier based fictitious domain methods  discussed in \cite{Glowinski1,Glowinski2,Glowinski3,Glowinski4, Glowinski5, Glowinski6}, which rely on $H^{1}$-- extensions, are not easy  to generalize to the slip situation. For incompressible viscous flow around obstacle with Navier slip boundary condition \cite{He-slip-Navier}, we have developed a new least--squares/fictitious domain method. The method in \cite{He-slip-Navier} is of the {\em virtual control} type (in the sense of J. L. Lions; see \cite{Lions}) and relies on a {\em least--squares} formulation making the problem solvable by a {\em conjugate gradient} algorithm operating in a well--chosen control space, which is a generation of least--squares/fictitious domain method for linear elliptic problem with Robin boundary  condition (see \cite{He1}) to Navier--Stokes problem with Navier slip boundary condition. Such least--squares idea has been sucessfully generated to Navier--Stokes--Cahn--Hilliard system in \cite{He-cahn}.

The main goal of the present article is to discuss the numerical simulation of fluid-particle interaction when Navier slip take effect at the interface fluid/particle. The particle velocity is updated together with the solution of Navier--Stokes equations for the stability. The position of particle is updated explicitly. The formulation of the problem is given in Section 2. In Section 3, we describe an operator--splitting scheme for the problem, which including the discription of  a least--squares/fictitious domain method for the sub-problem and a conjugate gradint solution of the   least--squares problem. 
The finite element implementation of the above methodology is discussed in Section 4.
Finally, we present in Section 5 the results of numerical experiments.

\section{\bf{Problem formulation}}
\begin{figure}[htp]
	\centering
	\includegraphics[height=2.8in]{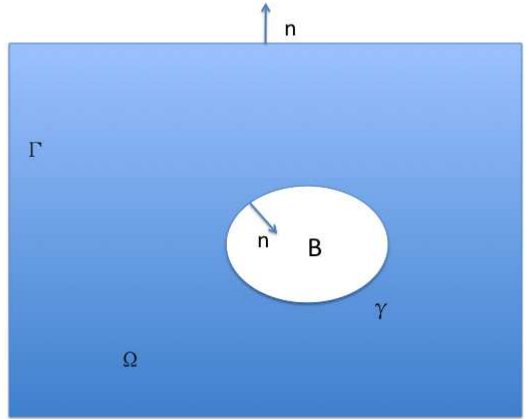}
	\caption{Problem geometry}
	\label{fig:1}     
\end{figure}

\setcounter{equation}{0}
Let $\Omega$ be a bounded domain of $\mathbb{R}^{d} (d = 2$ or 3), including the interior of the particle. We denote by $\Gamma$ the boundary of $\Omega$. The geometry of the problem has been visualized on Figure \ref{fig:1}. We suppose that $\Omega$ contains a rigid body $B$ whose boundary  is denoted by $\gamma$.  For simplicity, we assume that the fluid velocity satisfies a Dirichlet boundary condition on the outer boundary $\Gamma$. The $\rho_{s}$  and $\rho_{f}$ are particle density and fluid density, respectively. Compared to  previous work, the main novelty here is that we will assume that a slip boundary condition holds at the interface $\gamma$ between $\Omega$ and $B$. Assuming that the external force is $\mathbf{g}$, the motion of fluid--particle mixture is governed by the following equations:\\
\noindent
{\em Fluid motion}
\begin{align}
&  \rho_{f}\left[ \frac{\partial \mathbf{u}}{\partial t} + (\mathbf{u} \cdot \nabla) \mathbf{u}\right]  = \rho_{f} \mathbf{g} + \nabla \cdot \pmb{\sigma}, \quad  \mbox{in} \quad  \Omega\setminus \overline{B(t)}\label{eq:problemformulation1} \\
&    \nabla  \cdot \mathbf{u} = 0 \quad  \mbox{in} \quad \Omega\setminus \overline{B(t)}\label{eq:problemformulation2}\\
&  (\mathbf{u}-\mathbf{u}_{B})\cdot \mathbf{n} =0 \quad \mbox{on} \quad \gamma \times (0,T), \label{eq:problemformulation3}\\
&  (\pmb{\sigma}\mathbf{n}+\frac{\mu}{l_{s}} (\mathbf{u}-\mathbf{u}_{B}))\times \mathbf{n} = 0 \quad \mbox{on} \quad \gamma \times (0,T), \label{eq:problemformulation4}\\
&   \mathbf{u} = \mathbf{g}_{\Gamma} \quad  \mbox{on} \quad \Gamma \times (0,T), \label{eq:problemformulation5}\\
& \mathbf{u}_{B} = \mathbf{U} + \pmb{\omega}\times \mathbf{r} \ \mbox{on} \ \gamma,  \label{eq:problemformulation6}\\
&  \mathbf{u}(x,0) = \mathbf{u}_{0}, x\in \Omega \setminus \overline{ B(0)}, \quad \mbox{with} \quad \nabla  \cdot \mathbf{u_{0}} = 0.\label{eq:problemformulation7}
\end{align}
{\em Particle motion}
\begin{align}
&  \mathbf{M}\frac{d \mathbf{U}}{d t} = \mathbf{M}\mathbf{g} + \mathbf{F}, \label{eq:particle-motion1}\\
&  \frac{d\mathbf{I} \pmb{\omega}}{d t} = \mathbf{T}, \label{eq:particle-motion2} \\
&  \mathbf{U}|_{t=0} = \mathbf{U}_{0}, \label{eq:particle-motion3} \\
&  \pmb{\omega}|_{t=0} = \pmb{\omega}_{0} \label{eq:particle-motion4} \\
&  \frac{d \mathbf{X}}{d t} = \mathbf{U},\label{eq:particle-motion5}\\
&  \mathbf{X}|_{t = 0} = \mathbf{X}_{0}, \label{eq:particle-motion6}
\end{align}
where  $\mathbf{n}$ denotes the outward normal unit vector at $\partial (\Omega \backslash \overline{B})$, pointing outward to the fluid region,  $\mathbf{r} \triangleq \mathbf{x}-\mathbf{X}$; $\mathbf{u}$ and $\pmb{\sigma}$ are the fluid velocity and stress;   $\mathbf{M},\mathbf{I}, \mathbf{U}, \pmb{\omega}, \mathbf{X}$ are the mass, moment of inertia,translational velocity, angular velocity, center of mass of particle respectively;  and $ \mathbf{F}= -\int_{\partial B(t)}\pmb{\sigma}\mathbf{n} d \gamma$,  $\mathbf{T} = -\int_{\partial B(t)}\mathbf{r}\times\pmb{\sigma} \mathbf{n} d \gamma$ are the hydrodynamics force and torque on the  particle, 
and $l_{s} > 0$ is the slip length,  $\mathbf{u}_{B} $ is the velocity of the rigid body, $\mu$ is positive constant and denote  fluid viscosity.
Here $\pmb{\sigma} = 2 \mu \mathbf{D}(\mathbf{u})-p \mathbf{E}$ with $ \mathbf{D}(\mathbf{v})= \frac{1}{2} [\nabla \mathbf{v} +( \nabla \mathbf{v})^{t} ], \forall \ \mathbf{v}$, where $\mathbf{E}$ is identity matrix.

A classical variational formulation of the problem \eqref{eq:problemformulation1}--\eqref{eq:problemformulation7} is given by:
\noindent
Find $\{\mathbf{u}, p, \mathbf{U}, \pmb{\omega}\}$ such that, a.e. on $(0,T)$, \eqref{eq:formulate1}--\eqref{eq:formulate5} 
 are established.
\begin{align}
&  \rho_{f}\int_{\Omega\backslash \bar{B}}\left[\frac{\partial \mathbf{u}}{\partial t}+ (\mathbf{u} \cdot \nabla) \mathbf{u}\right]\cdot\mathbf{v}d x+2\mu \int_{\Omega\backslash\bar{B}}\mathbf{D}( \mathbf{u}) : \mathbf{D}  (\mathbf{v} )dx   \nonumber \\
&  + \int_{\Omega\backslash \bar{B}}\nabla p \cdot \mathbf{v} d x +  \frac{\mu}{l_{s}}\int_{\gamma} (\mathbf{u}-(\mathbf{U}+\pmb{\omega}\times \mathbf{r}))\cdot \mathbf{v} d \gamma = \rho_{f} \int_{\Omega\backslash \bar{B}}\mathbf{g}\cdot \mathbf{v} d x,  \label{eq:formulate1}\\
 &  \forall \mathbf{v} \in \mathbb{V}_{0}, \nonumber \\
&  \int_{\Omega\backslash \bar{B}} \nabla  \cdot \mathbf{u} \  q dx = 0, \  \forall q \in L^{2}(\Omega\backslash \bar{B}),\label{eq:formulate2}\\
&  \mathbf{u} = \mathbf{g}_{\Gamma} \  \mbox{on} \  \Gamma, \label{eq:formulate3}\\
&  \mathbf{u}(0) = \mathbf{u}_{0}, \label{eq:formulate4}\\
& \left(\mathbf{u}-(\mathbf{U}+\pmb{\omega}\times \mathbf{r})\right) \cdot \mathbf{n} = 0 \ \mbox{on} \ \gamma\times(0,T),\label{eq:formulate5}
\end{align}
where $ \mathbb{V}_{0} = \{ \mathbf{v}| \mathbf{v} \in (H^{1}(\Omega\backslash \bar{B}))^{d}, \ \mathbf{v} = \mathbf{0} \  \mbox{on} \ \Gamma, \  \mathbf{v} \cdot \mathbf{n} = 0 \ \mbox{on} \ \gamma  \}.$  

\section{\bf{An operator--splitting scheme for the time--discretization of problem \eqref{eq:formulate1}--\eqref{eq:formulate5} and \eqref{eq:particle-motion1}--\eqref{eq:particle-motion6}}}
We first  came up with a very simple scheme which decoupled the motion of the particle and the motion of fluid at each time step. That is, we use the method in \cite{He-slip-Navier} to solve  \eqref{eq:time-discretization2}--\eqref{eq:time-discretization5} firstly,  solve the pure advection problem \eqref{eq:time-discretization6}--\eqref{eq:time-discretization9} secondly,  update particle velocity and position explicitly. Unfortunately, this simple scheme is unstable. 
For $n \geq 0$,  $\mathbf{u}^{n}$ being known, 
$\Delta t > 0$ is a time--discretization step and $t^{n} = n \Delta t$,
\begin{eqnarray}
& & \rho_{f}\int_{\Omega\backslash \bar{B}}\left[\frac{\mathbf{u}^{n+\frac{1}{2}}-\mathbf{u}^{n}}{\Delta t}\cdot\mathbf{v}d x\right]+2\mu \int_{\Omega\backslash \bar{B}}\mathbf{D} (\mathbf{u}^{n+\frac{1}{2}}): \mathbf{D}( \mathbf{v} )dx + \int_{\Omega\backslash \bar{B}}\nabla p^{n+1} \cdot \mathbf{v} d x  \nonumber \\
& &  +  \frac{\mu}{l_{s}}\int_{\gamma} (\mathbf{u}^{n+\frac{1}{2}}-\mathbf{u}_{B})\cdot \mathbf{v} d \gamma = \rho_{f} \int_{\Omega\backslash \bar{B}}\mathbf{g}\cdot \mathbf{v} d x, \   \forall  \mathbf{v} \in \mathbb{V}_{0},  \label{eq:time-discretization2}  \\
& & \int_{\Omega\backslash \bar{B}} \nabla \cdot \mathbf{u}^{n+\frac{1}{2}} q dx = 0, \ \forall q \in L^{2}(\Omega\backslash \bar{B}),   \label{eq:time-discretization3}\\
& & (\mathbf{u}^{n+\frac{1}{2}}-\mathbf{u}_{B})\cdot \mathbf{n} = 0 \ \mbox{on} \ \gamma,   \label{eq:time-discretization4} \\
& & \mathbf{u}^{n+\frac{1}{2}} = \mathbf{g}_{\Gamma} \  \mbox{on} \ \Gamma,  \label{eq:time-discretization5}
\end{eqnarray}
where $\mathbf{u}_{B} = \mathbf{U} + \pmb{\omega} \times \mathbf{r}$ is obtained by \eqref{eq:particle-motion1} and \eqref{eq:particle-motion2}.
In order to solve problem \eqref{eq:formulate1}--\eqref{eq:formulate5} and \eqref{eq:particle-motion1}--\eqref{eq:particle-motion6} numerically, we advocate the
operator--splitting scheme where 
fluid velocity and particle velocty are updated implicitly together and particle position is updated explicitly.
 Then we solve the following pure advection problem
\begin{align}
&  \frac{\partial \mathbf{u}}{\partial t} + (\mathbf{u}^{n+\frac{1}{2}} \cdot \nabla )\mathbf{u} = \mathbf{0} \ \mbox{in} \ (\Omega \backslash \bar{B}) \times (t^{n}, t^{n+1}), \label{eq:time-discretization6}\\
&  \mathbf{u}(t^{n}) = \mathbf{u}^{n+\frac{1}{2}}, \label{eq:time-discretization7}\\
&  \mathbf{u} = \mathbf{g}_{\Gamma} \ \mbox{on} \  \mathbf{\Gamma}_{-}\times (t^{n}, t^{n+1}); \label{eq:time-discretization8}\\
&  \mathbf{u} = \mathbf{u}^{n+\frac{1}{2}} \ \mbox{on} \   \mathbf{\gamma}_{-}\times (t^{n}, t^{n+1}), \label{eq:time-discretization9}
\end{align}
set 
\begin{eqnarray}
\mathbf{u}^{n+1} = \mathbf{u}(t^{n+1}).\label{eq:time-discretization10}
\end{eqnarray}
where, $ \mathbf{\Gamma}_{-} = \{x| x\in \Gamma, \mathbf{g}_{\Gamma} \cdot \mathbf{n}  < 0\}$ and $ \mathbf{\gamma}_{-} = \{x| x\in \gamma, \mathbf{u}_{B} \cdot \mathbf{n}  < 0\}$. Since we consider low Reynolds number flows, the backward method of characteristics will be used to solve problems \eqref{eq:time-discretization6}--\eqref{eq:time-discretization9}. The implementation of such method has been detailed described in Chapter 6 in \cite{Handbook}. 
 Finally, the position of particle is updated by the approach in \cite{Glowinski6} for \eqref{eq:particle-motion5}--\eqref{eq:particle-motion6}. At each time step, keeping the distance constant between points $\mathbf{X}_{i}$ in particle  is important since we are dealing with rigid particle. Therefore, our paper will focus on how to solve \eqref{eq:time-discretization2}--\eqref{eq:time-discretization5} and updating $\mathbf{u}_{B}$.

\subsection{On the fictitious domain solution of sub-problem excluding advection terms}
In order to solve fluid velocity and particle velocity together, we advocate the following variant of \emph{virtual control/ fictitious domain method} discussed in  \cite{He1}  and \cite{He-slip-Navier} (below, $\mathbf{u}_{\ast}, \ p_{1\ast}$ are defined over the whole $\Omega$ and verify $\mathbf{u}_{\ast}|_{\Omega \backslash \bar{B}} = \mathbf{u}^{n}, \  p_{1\ast}|_{\Omega\backslash \bar{B}} = p^{n}, \mathbf{U}_{\ast} = \mathbf{U}^{n}, \omega_{\ast} = \omega^{n}$ ).  For simlicity, we consider two dimensional motion (the following numerical method can be generated to three dimension case straightforwardly), therefore angular velocity $\pmb{\omega}$ reduced to $\omega$ and angular velocity equation becomes 
\begin{align}
&\mathbf{I} \frac{d \omega}{ d t} = \mathbf{T}. \nonumber
\end{align}
Suppose there exists  $\mathbf{y} \in (L^{2}(B))^{2},
\mathbf{C}_{1}\in \mathbb{R}^{2}, C_{2}\in \mathbb{R} $ 
such that the following relations hold:
\begin{align}
&  \mathbf{u}_{1} \in (H^{1}(\Omega))^{2},  \mathbf{u}_{1} = \mathbf{g}_{\Gamma} \  \mbox{on} \  \Gamma,\  p_{1} \in H^{1}(\Omega), \nonumber \\
&  \rho_{f}\int_{\Omega}\frac{\mathbf{u}_{1}-\mathbf{u}_{\ast}}{\Delta t}\cdot\mathbf{v}d x-\rho_{f} \int_{B}\frac{\mathbf{y}}{\Delta t}\cdot \mathbf{v} dx 
 + 2 \mu \int_{\Omega}\mathbf{D} (\mathbf{u}_{1}): \mathbf{D} (\mathbf{v}) dx + \int_{\Omega}\nabla p_{1} \cdot \mathbf{v} d x \nonumber\\
&  = \rho_{f} \int_{\Omega\backslash \bar{B}}\mathbf{g}\cdot \mathbf{v} d x,  \ \forall \mathbf{v} \in (H_{0}^{1}(\Omega))^{2}, \label{eq:least1}\\
&  \int_{\Omega}\nabla \cdot \mathbf{u}_{1}  q dx = -\int_{\Omega} \mathbf{u}_{1} \cdot \nabla q d x+\int_{\Gamma}(\mathbf{g}_{\Gamma}\cdot \mathbf{n} )q d \Gamma =0, \  \forall q \in H^{1}(\Omega), \label{eq:least2}\\
&  \mathbf{u}_{2} \in (H^{1}(B))^{2}, \ p_{2} \in H^{1}(B), \  \left(\mathbf{u}_{2}-(\widetilde{\mathbf{U}}+\widetilde{\omega} \times \mathbf{r}) -\frac{1}{\Delta t}\left(\mathbf{C}_{1}+C_{2}\times\mathbf{r}\right)\right)\cdot \mathbf{n} = 0 \ \mbox{on} \ \gamma, \nonumber \\
&  \rho_{f} \int_{B}\frac{\mathbf{u}_{2}-\mathbf{u}_{\ast}-\mathbf{y}}{\Delta t}\cdot \mathbf{v} dx
+  2\mu \int_{B}\mathbf{D} (\mathbf{u}_{2}):\mathbf{D}(\mathbf{v} )dx+ \int_{B}\nabla p_{2} \cdot \mathbf{v} d x \nonumber\\
&= \frac{\mu}{l_{s}}\int_{\gamma}(\mathbf{u}_{1}-(\widetilde{\mathbf{U}}+\widetilde{\omega} \times \mathbf{r}))\cdot \mathbf{v} d \gamma -\frac{1}{\Delta t} \frac{\mu}{l_{s}}\int_{\gamma}(\mathbf{C}_{1}+C_{2}\times \mathbf{r}) \cdot \mathbf{v}d \gamma,   \label{eq:least3} \\
&  \forall \mathbf{v} \in \mathbb{V}_{0B} \triangleq \{ (H^{1}(B))^{2}, \mathbf{v} \cdot \mathbf{n} = 0 \ \mbox{on} \ \gamma \},  \nonumber\\
& \int_{B}\nabla \cdot \mathbf{u}_{2}  q   dx = -\int_{B} \mathbf{u}_{2} \cdot \nabla q d x+\int_{\gamma}(\widetilde{\mathbf{U}}+ \widetilde{\omega} \times \mathbf{r})\cdot \mathbf{n} q d \gamma  +\frac{1}{\Delta t}\int_{\gamma}(\mathbf{C}_{1}+C_{2} \times \mathbf{r})\cdot \mathbf{n} q d \gamma \nonumber \\
& = 0, \  \forall q \in H^{1}(B),  \label{eq:least4} \\
& \mathbf{U} = \widetilde{\mathbf{U}} + \frac{1}{\Delta t}\mathbf{C}_{1} , \label{eq:least5}\\
& \omega = \widetilde{\omega} + \frac{1}{\Delta t}C_{2}, \label{eq:least6}
\end{align}
where
\begin{align*} 
\widetilde{\mathbf{U}} = \mathbf{U}_{\ast} - \frac{\Delta t}{\mathbf{M}}\int_{\gamma}\left( -p_{1\ast}\mathbf{E}+2\mu\mathbf{D}(\mathbf{u}_{\ast})\right) \mathbf{n} d \gamma,
\end{align*}
 and 
 \begin{align*} 
 \widetilde{\omega} = \omega_{\ast}-\frac{\Delta t}{\mathbf{I}}\int_{\gamma}\mathbf{r}\times \left( -p_{1\ast}\mathbf{E}+2\mu\mathbf{D}(\mathbf{u}_{\ast})\right) \mathbf{n} d \gamma.
 \end{align*}
Both problems \eqref{eq:least1}--\eqref{eq:least2} and \eqref{eq:least3}--\eqref{eq:least4}  have a unique solution  in $(H^{1}(\Omega))^{2}\times H^{1}(\Omega)/\mathbb{R}$ and $(H^{1}(B))^{2}\times H^{1}(B)/\mathbb{R}$, respectively. We define $\mathbf{A}: (L^{2}(B))^{2}
 \times \mathbb{R}^{2} \times \mathbb{R}
 \rightarrow (H^{1}(B))^{2}$ by
\begin{eqnarray}
\mathbf{A}(\mathbf{y}, \mathbf{C}_{1}, C_{2}) =
\left(\mathbf{u}_{2}-\mathbf{u}_{1}\right)\mid_{B}. \label{eq:Boperator}
 \end{eqnarray}
 Operator  $\mathbf{A}$ is clearly {\em affine} and {\em continuous}. We observe that 
 if $\mathbf{y}^{\ast}, \mathbf{C}^{\ast}_{1}, C^{\ast}_{2}$
  verify 
  $\mathbf{A}(\mathbf{y},\mathbf{C}_{1}, C_{2}) = \mathbf{0}$,
  we then have $\mathbf{u}_{2} = \mathbf{u}_{1}$ on $B$ and it is easy to see that  $\mathbf{u}_{1}|_{\Omega \setminus \bar {B}}, p_{1}|_{\Omega \setminus \bar {B}}, \mathbf{U}$ and $\omega$ are what we need.  We will numerically discuss the solution of the following functional equation
    \begin{eqnarray}
    \mathbf{A}(\mathbf{y}^{\ast},\mathbf{C}^{\ast}_{1},C^{\ast}_{2}) = \mathbf{0}. \label{eq:functioneq}
    \end{eqnarray}
\begin{remark}
 Problem \eqref{eq:functioneq} can be viewed as an exact controllability problem in the sense of \cite{Lions}.
\end{remark}

\subsubsection{ A least--squares formulation}
In order to solve system \eqref{eq:functioneq}, we use the following least--squares approach:
\begin{eqnarray}
& & \mbox{Find} \  \mathbf{y}^{\ast}\in (L^{2}(B))^{2}, \mathbf{C}_{1}^{\ast}, C_{2}^{\ast} \ 
\mbox{such that}  \nonumber\\
& & J(\mathbf{y}^{\ast},\mathbf{C}^{\ast}_{1}, C^{\ast}_{2}  ) \leq
J(\mathbf{y},\mathbf{C}_{1}, C_{2}), \ \forall \mathbf{y}\in (L^{2}(B))^{2},  \mathbf{C}_{1}\in \mathbb{R}^{2}, \ C_{2}\in\mathbb{R},
 \label{eq:Jdefine1}
\end{eqnarray}
where 
\begin{eqnarray}
& & J(\mathbf{y},\mathbf{C}_{1}, C_{2} ) = 
\frac{1}{2}
\left[ \rho_{f}\int_{B}|\mathbf{u}_{2}-\mathbf{u}_{1}|^{2} d x +2 \mu \Delta t \int_{B}|\mathbf{D}(\mathbf{u}_{2}-\mathbf{u}_{1})|^{2}d x\right], \label{eq:Jdefine2}
\end{eqnarray}
where in \eqref{eq:Jdefine2}, $\mathbf{u}_{1}$ and $\mathbf{u}_{2}$ are obtained from $\mathbf{y}, \mathbf{C}_{1}$ and $C_{2}$, via the solutions of  \eqref{eq:least1}--\eqref{eq:least2} and  \eqref{eq:least3}--\eqref{eq:least4}. A natural
candidate for the solution of the minimization problem \eqref{eq:Jdefine1}--\eqref{eq:Jdefine2} is a conjugate gradient algorithm operating in the virtual control space $(L^{2}(B))^{2}\times \mathbb{R}^{2} \times \mathbb{R}$. 
Such a solution is characterized by 
\begin{eqnarray*}
\frac{\partial J}{\partial \mathbf{y}}(\mathbf{y}^{\ast},\mathbf{C}^{\ast}_{1},C^{\ast}_{2}) =
 \mathbf{0},  \  
\frac{\partial J}{\partial \mathbf{C}_{1}}(\mathbf{y}^{\ast},\mathbf{C}^{\ast}_{1},C^{\ast}_{2})= \mathbf{0}, \  \frac{\partial J}{\partial C_{2}}(\mathbf{y}^{\ast},\mathbf{C}^{\ast}_{1},C^{\ast}_{2})= 0.
\end{eqnarray*}

\subsubsection{On the computation of partial differentials of functional $J$}
We are going to address this most important issue using a \emph{perturbation  approach} (as done in \cite{He1}  for a close variant of problem \eqref{eq:Jdefine1}, and in \cite{Lions} for various linear and nonlinear control problems). Suppose that $\mathbf{y} \in (L^{2}(B))^{2},
\mathbf{C}_{1} \in \mathbb{R}^{2}$ and $C_{2} \in \mathbb{R}$, 
 a \emph{perturbation} $ \delta \mathbf{y} $ also belonging to $(L^{2}(B))^{2},
   \delta \mathbf{C}_{1}$ belonging to $\mathbb{R}^{2}$ and $ \delta C_{2}$  belonging to $\mathbb{R}$ are given,
 we have then 
\begin{align}
&  \delta J(\mathbf{y}, \mathbf{C}_{1}, C_{2}) = \rho_{f}\int_{B}
 (\mathbf{u}_{2}-\mathbf{u}_{1})\cdot \delta (\mathbf{u}_{2}-\mathbf{u}_{1}) d x  
  +  2\mu \Delta t \int_{B}\mathbf{D} (\mathbf{u}_{2}-\mathbf{u}_{1}): \mathbf{D} \delta( \mathbf{u}_{2}-\mathbf{u}_{1} )dx  \nonumber \\
& =  \rho_{f}\int_{B}(\mathbf{u}_{2}-\mathbf{u}_{1}) \delta \mathbf{u}_{2} dx  + 2\mu \Delta t \int_{B}\mathbf{D} (\mathbf{u}_{2}-\mathbf{u}_{1}): \mathbf{D}( \delta \mathbf{u}_{2} )dx\nonumber \\
&  +  \rho_{f}\int_{B}(\mathbf{u}_{1}-\mathbf{u}_{2}) \delta \mathbf{u}_{1} dx  + 2\mu \Delta t \int_{B}\mathbf{D} (\mathbf{u}_{1}-\mathbf{u}_{2}): \mathbf{D}( \delta \mathbf{u}_{1} )dx, \label{eq:Deferential}
\end{align}
with $\delta \mathbf{u}_{1}$ and $\delta \mathbf{u}_{2}$ verifying
\begin{align}
&  \delta \mathbf{u}_{1} \in (H^{1}_{0}(\Omega))^{2}, \ \delta p_{1} \in H^{1}(\Omega), \nonumber\\
&   \rho_{f}\int_{\Omega}\delta \mathbf{u}_{1}\cdot\mathbf{v}d x
 + 2 \mu \Delta t\int_{\Omega}\mathbf{D} (\delta \mathbf{u}_{1}): \mathbf{D} (\mathbf{v}) dx 
+ \Delta t \int_{\Omega}\nabla \delta p_{1} \cdot \mathbf{v} d x \nonumber \\
& =   \rho_{f} \int_{B}\delta \mathbf{y}\cdot \mathbf{v} d x, \  \forall \mathbf{v}\in (H^{1}_{0}(\Omega))^{2}, \label{eq:deltaU1} \\
&  \int_{\Omega}  \delta\mathbf{u}_{1} \cdot \nabla q d x = 0, \ \forall q \in H^{1}(\Omega), \label{eq:deltaU1p}\\
&  \delta \mathbf{u}_{2} \in (H^{1}(B))^{2}, \  \delta p_{2} \in H^{1}(B),\ \left(\delta \mathbf{u}_{2}-\frac{1}{\Delta t}\left(\delta \mathbf{C}_{1}+\delta C_{2}\times\mathbf{r}\right)\right)\cdot \mathbf{n} = 0 \ \mbox{on} \ \gamma, \nonumber\\
&   \rho_{f}\int_{B}\delta \mathbf{u}_{2}\cdot\mathbf{v}d x
 + 2 \mu \Delta t\int_{B}\mathbf{D} (\delta \mathbf{u}_{2}): \mathbf{D} (\mathbf{v}) dx 
+ \Delta t\int_{B}\nabla \delta p_{2} \cdot \mathbf{v} d x \nonumber \\
& =  \frac{\mu \Delta t}{l_{s}} \int_{\gamma} \delta \mathbf{u}_{1}\cdot \mathbf{v}d \gamma - \frac{\mu}{l_{s}}
\int_{\gamma} \left(\delta \mathbf{C}_{1}+ \delta C_{2}\times \mathbf{r}\right)\cdot \mathbf{v}d \gamma
+ \rho_{f} \int_{B}\delta \mathbf{y}\cdot \mathbf{v} d x, \  \forall \mathbf{v}\in \mathbb{V}_{0B}, \label{eq:deltaU2} \\
& \int_{B}  \delta\mathbf{u}_{2} \cdot \nabla q d x =  
 \frac{1}{\Delta t}\int_{\gamma} (\delta \mathbf{C}_{1}+\delta C_{2}\times \mathbf{r}) \cdot \mathbf{n} q d \gamma, 
\ \forall q \in H^{1}(B). \label{eq:deltaU2p}
\end{align} 
Let us define the function $\pmb{\pi}_{2}, \  p_{\pmb{\pi}_ {2}}$ as the solution of the following linear variational problem:
\begin{align}
&  \pmb{\pi}_{2} \in \mathbb{V}_{0B}, \  p_{\pmb{\pi}_ {2}} \in H^{1}(B), \nonumber \\
&  \rho_{f}\int_{B} \pmb{\pi}_{2}\cdot \mathbf{v} d x + 2 \mu \Delta t\int_{B} \mathbf{D}(\pmb{\pi}_{2}): \mathbf{D} (\mathbf{v})d x+ \Delta t\int_{B} \nabla
p_{\pmb{\pi}_ {2}} \cdot \mathbf{v} d x \nonumber \\
& = \rho_{f}\int_{B}(\mathbf{u}_{2}-\mathbf{u}_{1})\cdot \mathbf{v} d x+ 2\mu \Delta t \int_{B} \mathbf{D}(\mathbf{u}_{2}-\mathbf{u}_{1})\cdot \mathbf{D}(\mathbf{v})d x, \ \forall \mathbf{v} \in \mathbb{V}_{0B},  \label{eq:pi2}\\
& \int_{B} \pmb{\pi}_ {2} \cdot \nabla q d x = 0, \ \forall  q \in H^{1}(B).\label{eq:Ppi2}
\end{align}
Taking $\mathbf{v} = \delta \mathbf{u}_{2}- \frac{1}{\Delta t}\left(\delta \mathbf{C}_{1}+\delta C_{2}\times \mathbf{r}\right)$ 
in \eqref{eq:pi2} and $\mathbf{v} = \pmb{\pi}_ {2}$ in \eqref{eq:deltaU2}, using $\mathbf{D}(\delta \mathbf{C}_{1}+\delta C_{2}\times \mathbf{r}) = \mathbf{0}$, and 
combining with \eqref{eq:Deferential}, we obtain
\begin{align}
 &\delta J(\mathbf{y},\mathbf{C}_{1}, C_{2})   =  \rho_{f}\int_{B} \pmb{\pi}_{2}\cdot \left(\delta \mathbf{u}_{2} -\frac{1}{\Delta t}(\delta \mathbf{C}_{1}+\delta C_{2}\times \mathbf{r})\right) d x  \nonumber \\
&+ 2\mu\Delta t \int_{B} \mathbf{D}(\pmb{\pi}_{2}): \mathbf{D}\left( \delta \mathbf{u}_{2}-\frac{1}{\Delta t}(\delta \mathbf{C}_{1}+\delta C_{2}\times \mathbf{r})\right)d x   \nonumber\\
&+ \Delta t\int_{B} \nabla p_{\pmb{\pi}_ {2}} \cdot \left(\delta \mathbf{u}_{2}-\frac{1}{\Delta t}(\delta \mathbf{C}_{1}+\delta C_{2}\times \mathbf{r}) \right) dx
  + \frac{\rho_{f}}{\Delta t}\int_{B}(\mathbf{u}_{2}-\mathbf{u}_{1})\cdot (\delta \mathbf{C}_{1}+\delta C_{2}\times \mathbf{r}) d x \nonumber \\
 & + \rho_{f}\int_{B}(\mathbf{u}_{1}-\mathbf{u}_{2}) \cdot\delta \mathbf{u}_{1} dx  + 2\mu \Delta t \int_{B}\mathbf{D} (\mathbf{u}_{1}-\mathbf{u}_{2}): \mathbf{D}( \delta \mathbf{u}_{1} )dx \nonumber \\
& =   \rho_{f} \int_{B} \delta \mathbf{y} \cdot \pmb{\pi}_ {2} dx + \frac{\mu}{l_{s}} \int_{\gamma} \left(\Delta t\delta \mathbf{u}_{1}-(\delta \mathbf{C}_{1}+\delta C_{2}\times \mathbf{r})\right)\cdot \pmb{\pi}_{2} d \gamma- \Delta t\int_{B} \nabla \delta p_{2} \cdot \pmb{\pi}_{2} dx \nonumber \\
& + \Delta t\int_{B} \nabla p_{\pmb{\pi}_ {2}} \cdot \left(\delta \mathbf{u}_{2}- \frac{1}{\Delta t}(\delta \mathbf{C}_{1}+\delta C_{2}\times \mathbf{r})\right)dx + \frac{\rho_{f}}{\Delta t}\int_{B}(\mathbf{u}_{2}-\mathbf{u}_{1}-\pmb{\pi}_{2})\cdot(\delta \mathbf{C}_{1}+\delta C_{2}\times \mathbf{r})d x\nonumber \\
 &+  \rho_{f}\int_{B}(\mathbf{u}_{1}-\mathbf{u}_{2}) \cdot\delta \mathbf{u}_{1} dx  + 2\mu \Delta t \int_{B}\mathbf{D}  (\mathbf{u}_{1}-\mathbf{u}_{2}): \mathbf{D}( \delta \mathbf{u}_{1} )dx. \label{eq:deltaeq1}
\end{align}
If $q = \delta p_{2}$ in \eqref{eq:Ppi2} and $q =p_{\pmb{\pi}_ {2}}  $ in \eqref{eq:deltaU2p},  
 we have
\begin{align}
&\delta J(\mathbf{y},\mathbf{C}_{1},C_{2})  =    \rho_{f} \int_{B} \delta \mathbf{y} \cdot \pmb{\pi}_ {2} dx + \frac{\mu }{l_{s}} \int_{\gamma} \left(\Delta t\delta \mathbf{u}_{1}-(\delta \mathbf{C}_{1}+\delta C_{2}\times \mathbf{r})\right)\cdot \pmb{\pi}_{2} d \gamma \nonumber \\
&+ \frac{\rho_{f}}{\Delta t}\int_{B}(\mathbf{u}_{2}-\mathbf{u}_{1}-\pmb{\pi}_{2})\cdot \left(\delta \mathbf{C}_{1}+\delta C_{2}\times \mathbf{r})\right) dx \nonumber\\
&+ \rho_{f}\int_{B}(\mathbf{u}_{1}-\mathbf{u}_{2})\cdot \delta \mathbf{u}_{1} dx  + 2\mu \Delta t \int_{B}\mathbf{D} (\mathbf{u}_{1}-\mathbf{u}_{2}): \mathbf{D}( \delta \mathbf{u}_{1} )dx.\label{eq:deltaeq2}
\end{align}
We introduce now $\pmb{\pi}_{1}$ and $p_{\pmb{\pi}_ {1}}$ solution of the following linear variational problem:
\begin{align}
&  \pmb{\pi}_{1} \in (H^{1}_{0}(\Omega))^{2}, \  p_{\pmb{\pi}_ {1}} \in H^{1}(\Omega), \nonumber \\
 &\rho_{f}\int_{\Omega} \pmb{\pi}_{1}\cdot \mathbf{v} d x + 2 \mu \Delta t\int_{\Omega} \mathbf{D}(\pmb{\pi}_{1}): \mathbf{D} (\mathbf{v})d x+ \Delta t\int_{\Omega} \nabla
p_{\pmb{\pi}_ {1}} \cdot \mathbf{v} d x  \label{eq:pi1} \\
& =  \rho_{f}\int_{B}(\mathbf{u}_{1}-\mathbf{u}_{2})\cdot \mathbf{v} d x+ 2\mu \Delta t \int_{B} \mathbf{D}(\mathbf{u}_{1}-\mathbf{u}_{2})\cdot \mathbf{D}(\mathbf{v})d x+\frac{\mu \Delta t}{l_{s}}\int_{\gamma}\pmb{\pi}_{2}\cdot \mathbf{v} d \gamma, \nonumber \\
&   \forall \mathbf{v} \in (H^{1}_{0}(\Omega))^{2},  \nonumber\\
&  \int_{\Omega} \pmb{\pi}_ {1} \cdot \nabla q d x = 0, \ \forall  q \in H^{1}(\Omega).\label{eq:Ppi1}
\end{align}
Taking $\mathbf{v} = \delta \mathbf{u}_{1}$ in \eqref{eq:pi1} and combining with \eqref{eq:deltaeq2}, we obtain 
 \begin{align}
&\delta J(\mathbf{y},\mathbf{C}_{1},C_{2})    =   \rho_{f} \int_{B} \delta \mathbf{y} \cdot \pmb{\pi}_ {2} dx  + \rho_{f}\int_{\Omega} \pmb{\pi}_{1}\cdot \delta \mathbf{u}_{1} d x + 2 \mu \Delta t\int_{\Omega} \mathbf{D}(\pmb{\pi}_{1}): \mathbf{D} (\delta \mathbf{u}_{1})d x \nonumber \\
& + \Delta t\int_{\Omega} \nabla
p_{\pmb{\pi}_ {1}} \cdot \delta \mathbf{u}_{1}d x +\frac{\rho_{f}}{\Delta t}\int_{B}(\mathbf{u}_{2}-\mathbf{u}_{1}-\pmb{\pi}_{2})\cdot(\delta \mathbf{C}_{1}+\delta C_{2}\times \mathbf{r})d x- \frac{\mu }{l_{s}} \int_{\gamma}  (\delta \mathbf{C}_{1}+\delta C_{2}\times \mathbf{r})\cdot \pmb{\pi}_{2} d \gamma.
 \label{eq:deltaeq20} 
\end{align}
If $q = p_{\pmb{\pi}_ {1}}$ in \eqref{eq:deltaU1p}, we have
\begin{align}
 &\delta J(\mathbf{y},\mathbf{C}_{1},C_{2})    =   \rho_{f} \int_{B} \delta \mathbf{y} \cdot \pmb{\pi}_ {2} dx  + \rho_{f}\int_{\Omega} \pmb{\pi}_{1}\cdot \delta \mathbf{u}_{1} d x + 2 \mu \Delta t\int_{\Omega} \mathbf{D}(\pmb{\pi}_{1}): \mathbf{D} (\delta \mathbf{u}_{1})d x \nonumber\\
&+\frac{\rho_{f}}{\Delta t}\int_{B}(\mathbf{u}_{2}-\mathbf{u}_{1}-\pmb{\pi}_{2})\cdot(\delta \mathbf{C}_{1}+\delta C_{2}\times \mathbf{r})d x- \frac{\mu }{l_{s}} \int_{\gamma}  (\delta \mathbf{C}_{1}+\delta C_{2}\times \mathbf{r})\cdot \pmb{\pi}_{2} d \gamma. \label{eq:deltaeq23} 
\end{align}
Taking $\mathbf{v} = \pmb{\pi}_{1}$ in \eqref{eq:deltaU1} and combining with \eqref{eq:deltaeq23}, we finally obtain
\begin{align}
&\delta J(\mathbf{y},\mathbf{C}_{1},C_{2})    =   \rho_{f} \int_{B} \delta \mathbf{y} \cdot \pmb{\pi}_ {2} dx + \rho_{f}\int_{B} \delta \mathbf{y}\cdot \pmb{\pi}_{1} dx 
-\Delta t\int_{\Omega}\nabla \delta p_{1}\cdot \pmb{\pi}_{1} d x \nonumber\\
& +\frac{\rho_{f}}{\Delta t}\int_{B}(\mathbf{u}_{2}-\mathbf{u}_{1}-\pmb{\pi}_{2})\cdot(\delta \mathbf{C}_{1}+\delta C_{2}\times \mathbf{r})d x- \frac{\mu}{l_{s}} \int_{\gamma}  (\delta \mathbf{C}_{1}+\delta C_{2}\times \mathbf{r})\cdot \pmb{\pi}_{2} d \gamma. \label{eq:deltaeq33}
\end{align}
If $q = \delta p_{1}  $ in \eqref{eq:Ppi1}, we have
\begin{align}
&\delta J(\mathbf{y},\mathbf{C}_{1},C_{2})  
= \rho_{f}\int_{B} (\pmb{\pi}_{1}+\pmb{\pi}_{2})\cdot \delta \mathbf{y} dx 
 +\frac{\rho_{f}}{\Delta t}\int_{B}(\mathbf{u}_{2}-\mathbf{u}_{1}-\pmb{\pi}_{2})\cdot(\delta \mathbf{C}_{1}+\delta C_{2}\times \mathbf{r})d x \nonumber\\
& - \frac{\mu}{l_{s}}\int_{\gamma} (\delta \mathbf{C}_{1}+\delta C_{2}\times \mathbf{r}) \cdot \pmb{\pi}_{2}d \gamma. \label{eq:deltaeq4} 
\end{align}
Therefore, 
\begin{align}
 &\frac{\partial J}{\partial \mathbf{y}}(\mathbf{y},\mathbf{C}_{1},C_{2}) = \rho_{f}\left(\pmb{\pi}_{1}|_{B}+\pmb{\pi}_{2}\right), 
  \label{eq:deltaeq51} \\
&  \frac{\partial J}{\partial \mathbf{C}_{1}}(\mathbf{y},\mathbf{C}_{1},C_{2}) = \frac{\rho_{f}}{\Delta t}\int_{B}(\mathbf{u}_{2}-\mathbf{u}_{1}-\pmb{\pi}_{2})d x-\frac{\mu }{l_{s}} \int_{\gamma} \pmb{\pi}_{2} d \gamma, \label{eq:deltaeq52}  \\
& \frac{\partial J}{\partial C_{2}}(\mathbf{y},\mathbf{C}_{1},C_{2}) = \frac{\rho_{f}}{\Delta t}\int_{B}\left(-(\mathbf{u}_{2}-\mathbf{u}_{1}-\pmb{\pi}_{2})_{x}\mathbf{r}_{y}+(\mathbf{u}_{2}-\mathbf{u}_{1}-\pmb{\pi}_{2})_{y}\mathbf{r}_{x}\right) dx \nonumber\\
& -\frac{\mu}{l_{s}} \int_{\gamma}\left(-(\pmb{\pi}_{2})_{x}\mathbf{r}_{y}+(\pmb{\pi}_{2})_{y}\mathbf{r}_{x}\right) d\gamma, \label{eq:deltaeq53}
\end{align} 
where $(\cdot)_{x}$ and $(\cdot)_{y}$ are the x-component and y-component, respectively. 

\subsubsection{A conjugate gradient algorithm to the solution of the least--squares problem \eqref{eq:Jdefine1}--\eqref{eq:Jdefine2}}
Taking into account the way $\frac{\partial J}{\partial \mathbf{y}}(\mathbf{y}, \mathbf{C}_{1}, C_{2})$,
$ \  \frac{\partial J}{\partial \mathbf{C}_{1}}(\mathbf{y},\mathbf{C}_{1},C_{2})$ and $\frac{\partial J}{\partial C_{2}}(\mathbf{y},\mathbf{C}_{1},C_{2}) $,
 we are going to give a more practical formulation.   The algorithm reads as
\begin{align}
&  \mathbf{y}^{0},\ \mathbf{C}^{0}_{1}, \ C_{2}^{0} \ \mbox{are given in} \  (L^{2}(B))^{2} 
\times \mathbb{R}^{2}\times \mathbb{R}, \label{eq:algorithm1}
\end{align}
solve the following saddle point problems 
\begin{align}
&  \mathbf{u}^{0}_{1} \in (H^{1}(\Omega))^{2}, \ \mathbf{u}^{0}_{1}= \mathbf{g}_{\Gamma} \ \mbox{on} \ \Gamma, \ p^{0}_{1} \in H^{1}(\Omega), \nonumber\\
&  \rho_{f}\int_{\Omega}\mathbf{u}^{0}_{1}\cdot\mathbf{v} dx 
 + 2 \mu \Delta t\int_{\Omega}\mathbf{D} (\mathbf{u}^{0}_{1}): \mathbf{D} (\mathbf{v}) dx+ \Delta t\int_{\Omega}\nabla p^{0}_{1} \cdot \mathbf{v} d x \nonumber\\
& =  \rho_{f}\int_{B} \mathbf{y}^{0}\cdot \mathbf{v} d x+ \rho_{f}\Delta t\int_{\Omega\backslash \bar{B}}\mathbf{g}\cdot \mathbf{v} d x+\rho_{f} \int_{\Omega}\mathbf{u}_{\ast}\cdot \mathbf{v} d x,  \  \forall \mathbf{v} \in (H_{0}^{1}(\Omega))^{2}, \label{eq:algorithm2}\\
&  \int_{\Omega}\mathbf{u}^{0}_{1} \cdot \nabla  q dx = \int_{\Gamma}(\mathbf{g}_{\Gamma} \cdot \mathbf{n})  q d \Gamma , \  \forall q \in H^{1}(\Omega), \label{eq:algorithm3}\\
&  \mathbf{u}^{0}_{2} \in (H^{1}(B))^{2}, \  \left(\mathbf{u}^{0}_{2}-\left(\widetilde{\mathbf{U}}+\widetilde{\omega}\times \mathbf{r}\right)-\frac{1}{\Delta t}\left(\mathbf{C}^{0}_{1}+C^{0}_{2}\times\mathbf{r}\right)\right)\cdot \mathbf{n} =0 \ \mbox{on} \ \gamma, \ p^{0}_{2} \in H^{1}(B), \nonumber\\
&  \rho_{f}\int_{B}\mathbf{u}^{0}_{2}\cdot\mathbf{v} dx 
 + 2 \mu \Delta t\int_{B}\mathbf{D} (\mathbf{u}^{0}_{2}): \mathbf{D} (\mathbf{v}) dx+ \Delta t\int_{B}\nabla p^{0}_{2} \cdot \mathbf{v} d x 
 \nonumber \\
& =  \rho_{f}\int_{B} \mathbf{y}^{0}\cdot \mathbf{v} d x  +\rho_{f} \int_{B}\mathbf{u}_{\ast}\cdot \mathbf{v} d x \nonumber \\
&+\frac{\mu}{l_{s}} \int_{\gamma} \left(\Delta t\mathbf{u}^{0}_{1} -\Delta t\left(\widetilde{\mathbf{U}}+\widetilde{\omega}\times \mathbf{r}\right)-\left(\mathbf{C}^{0}_{1}+C^{0}_{2}\times\mathbf{r}\right)\right)\cdot \mathbf{v}  d \gamma,  \ \forall \mathbf{v} \in \mathbb{V}_{0B},  \label{eq:algorithm4}
\end{align}
\begin{align}
&  \int_{B}\mathbf{u}^{0}_{2} \cdot \nabla  q dx = \int_{\gamma}\left(\left(\widetilde{\mathbf{U}}+\widetilde{\omega}\times \mathbf{r}\right)+\frac{1}{\Delta t}\left(\mathbf{C}^{0}_{1}+C^{0}_{2}\times\mathbf{r}\right)\right)\cdot  \mathbf{n}q d \gamma , \  \forall q \in H^{1}(B), \label{eq:algorithm5}\\
&  \pmb{\pi}^{0}_{2} \in \mathbb{V}_{0B}, \   p^{0}_{\pmb{\pi}_{2}} \in H^{1}(B), \nonumber\\
&  \rho_{f}\int_{B}\pmb{\pi}^{0}_{2}\cdot\mathbf{v} dx 
 + 2 \mu \Delta t\int_{B}\mathbf{D} (\pmb{\pi}^{0}_{2}): \mathbf{D} (\mathbf{v}) dx+ \Delta t\int_{B}\nabla p^{0}_{\pmb{\pi}_{2}} \cdot \mathbf{v} d x  \label{eq:algorithm6} \\
& =  \rho_{f}\int_{B} (\mathbf{u}^{0}_{2}-\mathbf{u}^{0}_{1})\cdot \mathbf{v} d x+ 2\mu\Delta t \int_{B}\mathbf{D}(\mathbf{u}^{0}_{2}-\mathbf{u}^{0}_{1}): \mathbf{D}(\mathbf{v}),   \forall \mathbf{v} \in \mathbb{V}_{0B},\nonumber\\
&  \int_{B}\pmb{\pi}^{0}_{2} \cdot \nabla  q dx =0 , \  \forall q \in H^{1}(B), \label{eq:algorithm7}\\
&  \pmb{\pi}^{0}_{1} \in (H^{1}_{0}(\Omega))^{2}, \  p^{0}_{\pmb{\pi}_{1}} \in H^{1}(\Omega), \nonumber\\
&  \rho_{f}\int_{\Omega}\pmb{\pi}^{0}_{1}\cdot\mathbf{v} dx 
 + 2 \mu \Delta t\int_{\Omega}\mathbf{D} (\pmb{\pi}^{0}_{1}): \mathbf{D} (\mathbf{v}) dx+ \Delta t\int_{\Omega}\nabla p^{0}_{\pmb{\pi}_{1}} \cdot \mathbf{v} d x \nonumber\\
& =  \rho_{f}\int_{B} (\mathbf{u}^{0}_{1}-\mathbf{u}^{0}_{2})\cdot \mathbf{v} d x+ 2\mu\Delta t \int_{B}\mathbf{D}(\mathbf{u}^{0}_{1}-\mathbf{u}^{0}_{2}): \mathbf{D}(\mathbf{v})+\frac{\mu}{l_{s}}\Delta t \int_{\gamma} \pmb{\pi}^{0}_{2} \cdot \mathbf{v}  d \gamma, \nonumber\\
&  \forall \mathbf{v} \in (H^{1}_{0}(\Omega))^{2}, \label{eq:algorithm8}  \\
&  \int_{\Omega}\pmb{\pi}^{0}_{1} \cdot \nabla  q dx =0 , \  \forall q \in H^{1}(\Omega), \label{eq:algorithm9}
\end{align}
and set 
\begin{align*}
& \mathbf{g}^{0}_{1} = \rho_{f}(\pmb{\pi}^{0}_{1}|_{B}+\pmb{\pi}^{0}_{2}),\\
& \mathbf{w}^{0}_{1} = \mathbf{g}^{0}_{1},\\
& \mathbf{g}^{0}_{2} = \frac{\rho_{f}}{\Delta t}\int_{B}(\mathbf{u}^{0}_{2}-\mathbf{u}^{0}_{1}-\pmb{\pi}^{0}_{2}) d x- \frac{\mu }{l_{s}} \int_{\gamma} \pmb{\pi}^{0}_{2}d \gamma,\\
& \mathbf{w}^{0}_{2} = \mathbf{g}^{0}_{2},\\
& g^{0}_{3} = \frac{\rho_{f}}{\Delta t}\int_{B}\left(-(\mathbf{u}^{0}_{2}-\mathbf{u}^{0}_{1}-\pmb{\pi}^{0}_{2})_{x}\mathbf{r}_{y}+(\mathbf{u}^{0}_{2}-\mathbf{u}^{0}_{1}-\pmb{\pi}^{0}_{2})_{y}\mathbf{r}_{x}\right) dx \nonumber\\
&-\frac{\mu}{l_{s}} \int_{\gamma}\left(-(\pmb{\pi}^{0}_{2})_{x}\mathbf{r}_{y}+(\pmb{\pi}^{0}_{2})_{y}\mathbf{r}_{x}\right) d\gamma, \\
& w^{0}_{3} = g^{0}_{3}.
\end{align*}
For $k \geq 0, \mathbf{y}^{k}, \mathbf{C}^{k}_{1}, C^{k}_{2}, \mathbf{g}^{k}_{1}, \mathbf{g}^{k}_{2}, g^{k}_{3}$ 
and $\mathbf{w}^{k}_{1}, \mathbf{w}^{k}_{2}, w^{k}_{3}$
 being known, we compute 
$\mathbf{y}^{k+1}, \mathbf{C}^{k+1}_{1}, C^{k+1}_{2}$, $\mathbf{g}^{k+1}_{1}, \mathbf{g}^{k+1}_{2}, g_{3}^{k+1}$ and $\mathbf{w}_{1}^{k+1}, \mathbf{w}_{2}^{k+1}, w_{3}^{k+1}$ as follows: \\
\noindent
solve the saddle point problems
\begin{align}
&  \overline{\mathbf{u}}^{k}_{1} \in (H^{1}_{0}(\Omega))^{2},  \ \overline{p}^{k}_{1} \in H^{1}(\Omega), \nonumber\\
&  \rho_{f}\int_{\Omega}\overline{\mathbf{u}}^{k}_{1}\cdot\mathbf{v} dx 
 + 2 \mu \Delta t\int_{\Omega}\mathbf{D} (\overline{\mathbf{u}}^{k}_{1}): \mathbf{D} (\mathbf{v}) dx+ \Delta t \int_{\Omega}\nabla \overline{p}^{k}_{1} \cdot \mathbf{v} d x \nonumber\\
& =  \rho_{f}\int_{B} \mathbf{w}^{k}_{1}\cdot \mathbf{v} d x,  \ \forall \mathbf{v} \in (H_{0}^{1}(\Omega))^{2}, \label{eq:algorithm10}\\
&  \int_{\Omega}\overline{\mathbf{u}}^{k}_{1} \cdot \nabla  q dx = 0 , \  \forall q \in H^{1}(\Omega), \label{eq:algorithm11}\\
&  \overline{\mathbf{u}}^{k}_{2} \in \mathbb{V}_{0B},  \  \overline{p}^{k}_{2} \in H^{1}(B), \nonumber\\
&  \rho_{f}\int_{B} \overline{\mathbf{u}}^{k}_{2}\cdot\mathbf{v} dx 
 + 2 \mu \Delta t\int_{B}\mathbf{D} (\overline{\mathbf{u}}^{k}_{2}): \mathbf{D} (\mathbf{v}) dx+ \Delta t\int_{B}\nabla \overline{p}^{k}_{2} \cdot \mathbf{v} d x  \nonumber \\
& =  \rho_{f}\int_{B} \mathbf{w}^{k}\cdot \mathbf{v} d x 
+\frac{\mu}{l_{s}} \int_{\gamma}\left(\Delta t\overline{\mathbf{u}}^{k}_{1}-\left(\mathbf{w}^{k}_{2}+w^{k}_{3}\times \mathbf{r}\right) \right)\cdot \mathbf{v}  d \gamma,  \ \forall \mathbf{v} \in \mathbb{V}_{0B}, \label{eq:algorithm12}\\
&  \int_{B}\overline{\mathbf{u}}^{k}_{2} \cdot \nabla  q dx =0, \  \forall q \in H^{1}(B), \label{eq:algorithm13}\\
&  \overline{\pmb{\pi}}^{k}_{2} \in \mathbb{V}_{0B}, \  \overline{p}^{k}_{\pmb{\pi}_{2}} \in H^{1}(B), \nonumber\\
&  \rho_{f}\int_{B}\overline{\pmb{\pi}}^{k}_{2}\cdot\mathbf{v} dx 
 + 2 \mu \Delta t\int_{B}\mathbf{D} (\overline{\pmb{\pi}}^{k}_{2}): \mathbf{D} (\mathbf{v}) dx+ \Delta t\int_{B}\nabla \overline{p}^{k}_{\pmb{\pi}_{2}} \cdot \mathbf{v} d x   \label{eq:algorithm14} \\
& =  \rho_{f}\int_{B} (\overline{\mathbf{u}}^{k}_{2}-\overline{\mathbf{u}}^{k}_{1})\cdot \mathbf{v} d x+ 2\mu\Delta t \int_{B}\mathbf{D}(\overline{\mathbf{u}}^{k}_{2}-\overline{\mathbf{u}}^{k}_{1}): \mathbf{D}(\mathbf{v}),   \forall \mathbf{v} \in \mathbb{V}_{0B}, \nonumber\\
&  \int_{B}\overline{\pmb{\pi}}^{k}_{2} \cdot \nabla  q dx
= 0,   \forall q \in H^{1}(B), \label{eq:algorithm15}\\ 
& \overline{\pmb{\pi}}^{k}_{1} \in (H^{1}_{0}(\Omega))^{2}, \  \overline{p}^{k}_{\pmb{\pi}_{1}} \in H^{1}(\Omega), \nonumber\\
& \rho_{f}\int_{\Omega}\overline{\pmb{\pi}}^{k}_{1}\cdot\mathbf{v} dx 
 + 2 \mu \Delta t\int_{\Omega}\mathbf{D} (\overline{\pmb{\pi}}^{k}_{1}): \mathbf{D} (\mathbf{v}) dx+ \Delta t \int_{\Omega}\nabla \overline{p}^{k}_{\pmb{\pi}_{1}} \cdot \mathbf{v} d x \nonumber\\
& = \rho_{f} \int_{B} (\overline{\mathbf{u}}^{k}_{1}-\overline{\mathbf{u}}^{k}_{2}) \cdot \mathbf{v} d x+ 2\mu\Delta t \int_{B}\mathbf{D}(\overline{\mathbf{u}}^{k}_{1}-\overline{\mathbf{u}}^{k}_{2}): \mathbf{D}(\mathbf{v})+\frac{\mu}{l_{s}}\Delta t \int_{\gamma} \overline{\pmb{\pi}}^{k}_{2} \cdot \mathbf{v}  d \gamma,  \nonumber\\
&  \forall \mathbf{v} \in (H^{1}_{0}(\Omega))^{2},  \label{eq:algorithm16} \\
&  \int_{\Omega} \overline{\pmb{\pi}}^{k}_{1} \cdot \nabla  q dx =0 , \  \forall q \in H^{1}(\Omega), \label{eq:algorithm17}
\end{align}
and set 
\begin{align}
& \overline{\mathbf{g}}^{k}_{1} =  \rho_{f}(\overline{\pmb{\pi}}^{k}_{1}|_{B}+\bar{\pmb{\pi}}^{k}_{2}),  \label{eq:algorithm18}\\
& \overline{\mathbf{g}}^{k}_{2} = \frac{\rho_{f}}{\Delta t}\int_{B}(\overline{\mathbf{u}}^{k}_{2}-\overline{\mathbf{u}}^{k}_{1}-\overline{\pmb{\pi}}^{k}_{2}) d x- \frac{\mu}{l_{s}} \int_{\gamma} \overline{\pmb{\pi}}^{k}_{2} d \gamma, \label{eq:algorithm181}\\
& \overline{g}^{k}_{3} = \frac{\rho_{f}}{\Delta t}\int_{B}\left(-\left(\overline{\mathbf{u}}^{k}_{2}-\overline{\mathbf{u}}^{k}_{1}-\overline{\pmb{\pi}}^{k}_{2}\right)_{x}\mathbf{r}_{y}
+(\overline{\mathbf{u}}^{k}_{2}-\overline{\mathbf{u}}^{k}_{1}-\overline{\pmb{\pi}}^{k}_{2})_{y}\mathbf{r}_{x}\right) dx \nonumber\\
& -\frac{\mu}{l_{s}} \int_{\gamma}\left(-(\overline{\pmb{\pi}}^{k}_{2})_{x}\mathbf{r}_{y}+(\overline{\pmb{\pi}}^{k}_{2})_{y}\mathbf{r}_{x}\right) d\gamma. \label{eq:algorithm182}
\end{align}
Compute
\begin{align}
&\lambda_{k}  =  \frac{\int_{B}|\mathbf{g}^{k}_{1}|^{2}d x + |\mathbf{g}^{k}_{2}|^{2}+ |g^{k}_{3}|^{2}}{\int_{B} \overline{\mathbf{g}}^{k}_{1} \cdot \mathbf{w}^{k}_{1} d x + \overline{\mathbf{g}}^{k}_{2}\cdot \mathbf{w}^{k}_{2} + \overline{g}^{k}_{3}w^{k}_{3}},\label{eq:algorithm19}\\
&\mathbf{y}^{k+1}  =  \mathbf{y}^{k}-\lambda_{k} \mathbf{w}^{k}_{1}, \label{eq:algorithm20} \\
& \mathbf{C}^{k+1}_{1} = \mathbf{C}^{k}_{1}-\lambda_{k} \mathbf{w}^{k}_{2}, \label{eq:algorithm21}\\
& C^{k+1}_{2} = C^{k}_{2}-\lambda_{k}w^{k}_{3}, \label{eq:algorithm22}\\
&\mathbf{g}^{k+1}_{1}  =  \mathbf{g}^{k}_{1}-\lambda_{k} \overline{\mathbf{g}}^{k}_{1}, \label{eq:algorithm23}\\
& \mathbf{g}^{k+1}_{2} = \mathbf{g}^{k}_{2}-\lambda_{k}\overline{\mathbf{g}}^{k}_{2}, \label{eq:algorithm24}\\
& g^{k+1}_{3} = g^{k}_{3}-\lambda_{k}\overline{g}^{k}_{3}. \label{eq:algorithm25}
\end{align}
If  $\frac{\int_{B}|\mathbf{g}^{k+1}_{1}|^{2}d x + |\mathbf{g}^{k+1}_{2}|^{2}+ |g^{k+1}_{3}|^{2}}{\int_{B}|\mathbf{g}^{0}_{1}|^{2}d x+ |\mathbf{g}^{0}_{2}|^{2}+ |g^{0}_{3}|^{2}} \leq \ \mbox{tol}$,  we take $\mathbf{y} = \mathbf{y}^{k+1}, \  \mathbf{C}_{1} = \mathbf{C}^{k+1}_{1}, $ and $ C_{2} = C_{2}^{k+1}$,  
which is substituted into \eqref{eq:least1}--\eqref{eq:least2} and \eqref{eq:least5}--\eqref{eq:least6}, we get 
\begin{align*}
& \mathbf{u} = \mathbf{u}_{1}|_{\Omega\backslash \overline{B}},\\
& \mathbf{U}  =  \widetilde{\mathbf{U}} + \frac{1}{\Delta t}\mathbf{C}_{1}, \\
& \omega = \widetilde{\omega} + \frac{1}{\Delta t}C_{2}.
\end{align*}
Otherwise, compute
\begin{align}
& \gamma_{k}  =  \frac{\int_{B}|\mathbf{g}^{k+1}_{1}|^{2}d x +|\mathbf{g}^{k+1}_{2}|^{2}+ |g^{k+1}_{3}|^{2} }{\int_{B}|\mathbf{g}^{k}_{1}|^{2}d x +|\mathbf{g}^{k}_{2}|^{2}+ |g^{k}_{3}|^{2}}, \label{eq:algorithm26}
\end{align}
and set
\begin{align}
& \mathbf{w}^{k+1}_{1} =  \mathbf{g}^{k+1}_{1} + \gamma_{k}\mathbf{w}^{k}_{1}, \label{eq:algorithm27} \\
& \mathbf{w}^{k+1}_{2} =  \mathbf{g}^{k+1}_{2} + \gamma_{k}\mathbf{w}^{k}_{2}, \label{eq:algorithm28} \\
& w^{k+1}_{3} = g^{k+1}_{3} + \gamma_{k} w^{k}_{3}. \label{eq:algorithm29}
\end{align}
Do $n+1 \rightarrow n$ and return to \eqref{eq:algorithm10}. 

By inspection of \eqref{eq:algorithm2}--\eqref{eq:algorithm3},  \eqref{eq:algorithm8}--\eqref{eq:algorithm9}, \eqref{eq:algorithm10}--\eqref{eq:algorithm11} and \eqref{eq:algorithm16}--\eqref{eq:algorithm17},  it is clear that 
the pressure term can be treated as a corresponding Lagrange multiplier in the space $H_{0}$, which is defined as $H_{0} = \{ q | q\in H^{1}(\Omega), \ \int_{\Omega} q d x = 0\}$. Similar arguments can be applied to \eqref{eq:algorithm4}--\eqref{eq:algorithm5},  \eqref{eq:algorithm6}--\eqref{eq:algorithm7}, \eqref{eq:algorithm12}--\eqref{eq:algorithm13} and \eqref{eq:algorithm14}--\eqref{eq:algorithm15}.  Therefore, a preconditioned conjugate gradient algorithm  is easy to implement and seems to have excellent convergence properties. 
\begin{remark}
 If $\mathbf{y}^{0}$, $\mathbf{C}^{0}_{1}$, $C^{0}_{2}$ are close to $\mathbf{y}$, $\mathbf{C}_{1}$, $C_{2}$  (the stopping criterion we used for algorithm \eqref{eq:algorithm1}--\eqref{eq:algorithm29}) may lead to more iterations than necessary. A more realistic stopping test is given by
 \begin{eqnarray*}
 \frac{\int_{B} |\mathbf{g}^{n+1}_{1}|^{2}d x +  |\mathbf{g}^{k+1}_{2}|^{2}+ |g^{k+1}_{3}|^{2}}{\max\left(\int_{B} |\mathbf{g}^{0}_{1}|^{2}d x +  |\mathbf{g}^{0}_{2}|^{2}+ |g^{0}_{3}|^{2}, \int_{B} |\mathbf{y}^{n+1}|^{2} d x + |\mathbf{C}^{k+1}_{1}|^{2} + |C^{k+1}_{2}|^{2} \right)} \leq \ \mbox{tol}.
 \end{eqnarray*}
 Other stopping criteria can be used.
 \end{remark}

\section{Finite element realization of the operator-splitting scheme \eqref{eq:formulate1}--\eqref{eq:formulate5} and \eqref{eq:particle-motion1}--\eqref{eq:particle-motion6}}
\setcounter{equation}{0}
\subsection{Generalities}
We describe in this section  the {\em finite element} implementation of the least--squares/ fictitious domain methodology discussed in above sections. 
We will assume that $\overline{B} \subset \Omega \subset \mathbb{R}^{2}$ and that $\Omega$ is convex and/or has a smooth boundary; similarly, we assume that $\gamma$ is smooth. Our approximation of choice will be the  \emph{Bercovier--Pironneau} one, discussed with details in Glowinski\cite{Handbook}(see also the references therein). For simplicity we still denote by $\Omega$ and $B$ the polygonal approximations of the above domains. From the triangulations $\mathcal{T}_{h_{1}}$ of  $\Omega$ and $\mathcal{T}_{h_{2}}$ of  $B$ we define the following finite dimensional spaces.   The pressure space $P_{h_{1}}$ and $P_{Bh_{2}}$  are
\begin{eqnarray}
& & P_{h_{1}} = \{q | q \in C^{0}(\overline{\Omega}),  q|_{T} \in P_{1}, \forall T \in  \mathcal{T}_{h_{1}}\}, \label{eq:F-EGp1} \\
& & P_{Bh_{2}} = \{q | q \in C^{0}(\overline{B}),  q|_{T} \in P_{1}, \forall T \in  \mathcal{T}_{h_{2}}\}, \label{eq:F-EGp2}\\
& & P_{0h_{1}} = \{q | q \in P_{h_{1}}, \  \int_{\Omega} q  dx  = 0\}, \label{eq:F-EGp3}\\
& & P_{0Bh_{2}} = \{q |  q \in P_{Bh_{2}}, \  \int_{B} q  dx  = 0\}, \label{eq:F-EGp4}
\end{eqnarray}
where $P_{1}$ being the space of the polynomials of two variables of degree $\leq 1$  and $h_{1}$  (resp., $h_{2}$) the length of the largest edge(s)  of the finite element triangulation $\mathcal{T}_{h_{1}}$  (resp., $\mathcal{T}_{h_{2}}$).  Next, we divide each triangle of $ \mathcal{T}_{h_{1}}$(resp. $ \mathcal{T}_{h_{2}}$) in four sub-triangles, by joining the mid-points and denote by $ \mathcal{T}_{h_{1}/2}$(resp. $ \mathcal{T}_{h_{2}/2}$) the resulting triangulation. From $ \mathcal{T}_{h_{1}/2}$(resp. $ \mathcal{T}_{h_{2}/2}$), we define the velocity spaces 
\begin{eqnarray}
& & \mathbb{V}_{h_{1}} = \{ \mathbf{v} | \mathbf{v} \in C^{0}(\overline{\Omega}), \mathbf{v}|_{T}\in (P_{1})^{2}, \forall T \in  \mathcal{T}_{h_{1}/2}\}, \label{eq:F-EG1} \\
& & \mathbb{V}_{0h_{1}} = \{ \mathbf{v} | \mathbf{v}\in  \mathbb{V}_{h_{1}},  \mathbf{v} = \mathbf{0} \  \mbox{on} \ \Gamma \},\label{eq:F-EG2}
\end{eqnarray}
and
\begin{eqnarray}
& & \mathbb{V}_{Bh_{2}} = \{ \mathbf{v} | \mathbf{v} \in C^{0}(\overline{B}), \mathbf{v}|_{T} \in (P_{1})^{2}, \forall T \in  \mathcal{T}_{h_{2}/2}\}. \label{eq:F-EG3}
\end{eqnarray}
 The finite dimensional spaces $\mathbb{V}_{h_{1}}, \mathbb{V}_{0h_{1}}$ and $\mathbb{V}_{Bh_{2}}$ are finite dimensional approximations to $(H^{1}(\Omega))^{2}, (H_{0}^{1}(\Omega))^{2}$  and $(H^{1}(B))^{2}$, respectively. Similarly, we will use $\mathbb{V}_{Bh_{2}}$  to approximate the control space $(L^{2}(B))^{2}$.  We use $\Sigma_{h_{1}}$ to denote the set of the vertices of $ \mathcal{T}_{h_{1}/2}$ and $\Sigma_{h_{2}}$ to denote the set of the vertices of $ \mathcal{T}_{h_{2}/2}$.

\subsection{Finite element approximation of the least--squares problem \eqref{eq:Jdefine1}--\eqref{eq:Jdefine2}}
To approximate the least--squares problem \eqref{eq:Jdefine1}--\eqref{eq:Jdefine2}, we suggest
\begin{align}
&  \mathbf{y}^{\ast}_{h} \in \mathbb{V}_{Bh_{2}},  \  \mathbf{C}^{\ast}_{1} \in \mathbb{R}^{2}, \ C^{\ast}_{2} \in \mathbb{R}, \nonumber\\
&  J_{h}(\mathbf{y}^{\ast}_{h}, \mathbf{C}^{\ast}_{1}, C^{\ast}_{2}) \leq J_{h}(\mathbf{y}_{h}, \mathbf{C}_{1}, C_{2}), \ \forall \mathbf{y}_{h}\in \mathbb{V}_{Bh_{2}}, \  \mathbf{C}_{1} \in \mathbb{R}^{2}, \   C_{2} \in \mathbb{R} \label{eq:JdefineDis1}
\end{align}
where
\begin{eqnarray}
J_{h}(\mathbf{y}_{h}, \mathbf{C}_{1}, C_{2}) = \frac{1}{2}\left[ \rho_{f}\int_{B}|\mathbf{u}_{2}-\Pi_{2}\mathbf{u}_{1}|^{2} d x +2 \mu \Delta t \int_{B}|\mathbf{D}(\mathbf{u}_{2}-\Pi_{2}\mathbf{u}_{1})|^{2}d x\right]. \label{eq:JdefineDis2}
\end{eqnarray}
In the above, $\mathbf{u}_{1}$ is solution to the following fully discrete approximate saddle point problem:
\begin{align}
&  \mathbf{u}_{1} \in \mathbb{V}_{h_{1}}, \ \mathbf{u}_{1}= \mathbf{g}_{\Gamma} \ \mbox{on} \ \Gamma, \ p_{1} \in P_{h_{1}}, \nonumber\\
&  \rho_{f}\int_{\Omega}\mathbf{u}_{1}\cdot\mathbf{v} dx 
 + 2 \mu \Delta t\int_{\Omega}\mathbf{D} (\mathbf{u}_{1}): \mathbf{D} (\mathbf{v}) dx+ \Delta t\int_{\Omega}\nabla p_{1} \cdot \mathbf{v} d x \nonumber\\
& =  \rho_{f}\int_{B} \mathbf{y}_{h}\cdot \Pi_{2}\mathbf{v} d x+\rho_{f} \int_{\Omega}\mathbf{u}_{\ast}\cdot \mathbf{v} d x+ \rho_{f} \Delta t\int_{\Omega\backslash \bar{B}}\mathbf{g}\cdot \mathbf{v} d x, \ \forall \mathbf{v} \in \mathbb{V}_{0h_{1}}, \label{eq:dis-saddle1}\\
&  \int_{\Omega}\mathbf{u}_{1} \cdot \nabla  q dx = \int_{\Gamma}(\mathbf{g}_{\Gamma} \cdot \mathbf{n})  q d \Gamma , \  \forall q \in P_{h_{1}}, \label{eq:dis-saddle2}
\end{align}
where  $\Pi_{2}: (C^{0}(\overline{\Omega}))^{2} \rightarrow \mathbb{V}_{h_{2}}$ is the interpolation operator  defined as follows
\begin{align}
 \Pi_{2}\mathbf{v}=\sum_{i = 1}^{N_{h_{2}}} \mathbf{v}(Y_{i}) w_{2i}, \  \forall \  \mathbf{v} \in (C^{0}(\overline{\Omega}))^{2}, \label{eq:interpolation}
\end{align}
where $\{ Y_{i}\}_{i = 1}^{N_{h_{2}}}$ being the set of the vertices of $\mathcal{T}_{h_{2}/2}$ and $w_{2i}$ the $P_{1}$- shape function in $\mathbb{V}_{h_{2}}$ associated with the vertex $Y_{i}$. Return to \eqref{eq:JdefineDis2}, the function $\mathbf{u}_{2}$ is the solution of the following saddle point problem
\begin{align}
&  \mathbf{u}_{2} \in \mathbb{V}_{Bh_{2}}, \  \left(\mathbf{u}_{2}-(\widetilde{\mathbf{U}}+\widetilde{\omega} \times \mathbf{r}) -\frac{1}{\Delta t}\left(\mathbf{C}_{1}+C_{2}\times\mathbf{r}\right)\right)\cdot \mathbf{n} =0 \ \mbox{on} \ \gamma, \ p_{2} \in P_{B_{h_{2}}}, \nonumber\\
&  \rho_{f}\int_{B}\mathbf{u}_{2}\cdot\mathbf{v} dx 
 + 2 \mu \Delta t\int_{B}\mathbf{D} (\mathbf{u}_{2}): \mathbf{D} (\mathbf{v}) dx+ \Delta t\int_{B}\nabla p_{2} \cdot \mathbf{v} d x \label{eq:dis-saddle3} \\
& =  \rho_{f}\int_{B} \mathbf{y}_{h}\cdot \mathbf{v} d x+\frac{\mu\Delta t}{l_{s}}\int_{\gamma}(\Pi_{2}\mathbf{u}_{1}-(\widetilde{\mathbf{U}}+\widetilde{\omega} \times \mathbf{r}))\cdot \mathbf{v} d \gamma - \frac{\mu }{l_{s}}\int_{\gamma}(\mathbf{C}_{1}+C_{2}\times \mathbf{r}) \cdot \mathbf{v}d \gamma \nonumber\\
& +\rho_{f} \int_{B}\Pi_{2}\mathbf{u}_{\ast}\cdot \mathbf{v} d x,  \forall \mathbf{v} \in \mathbb{V}_{0Bh_{2}},  \nonumber\\
&  \int_{B}\mathbf{u}_{2} \cdot \nabla  q dx = \int_{\gamma}\left(\left(\widetilde{\mathbf{U}}+\widetilde{\omega}\times \mathbf{r}\right)+\frac{1}{\Delta t}\left(\mathbf{C}_{1}+C_{2}\times\mathbf{r}\right)\right)\cdot  q d \gamma \  \forall q \in P_{Bh_{2}}, \label{eq:dis-saddle4}\\
& \mathbf{U} = \widetilde{\mathbf{U}} + \frac{1}{\Delta t}\mathbf{C}_{1}, \label{eq:disU}\\
& \omega = \widetilde{\omega} + \frac{1}{\Delta t}C_{2}, \label{eq:disomega}
\end{align}
where $\mathbb{V}_{0Bh_{2}} \triangleq \{\mathbf{v}| \mathbf{v} \in \mathbb{V}_{Bh_{2}},  \ \mathbf{v}\cdot \mathbf{n} = 0 \ \mbox{on} \ \gamma \}, \widetilde{\mathbf{U}} = \mathbf{U}_{\ast} - \frac{\Delta t}{\mathbf{M}}\int_{\gamma}\left( -\Pi_{2}p_{1\ast}\mathbf{E}+2\mu\mathbf{D}(\Pi_{2}\mathbf{u}_{\ast})\right) d \gamma, $  and
$\widetilde{\omega} = \omega_{\ast} - \frac{\Delta t}{\mathbf{I}}\int_{\gamma}\mathbf{r}\times\left( -\Pi_{2}p_{1\ast}\mathbf{E}+2\mu\mathbf{D}(\Pi_{2}\mathbf{u}_{\ast})\right) d \gamma$.

With the introduce of interpolation operator $\Pi_{2}$, a finite element discretization of the conjugate gradient algorithms \eqref{eq:algorithm1}--\eqref{eq:algorithm29} is easy to realize.  As mentioned above, the finite element discretization of \eqref{eq:algorithm2}--\eqref{eq:algorithm3},  \eqref{eq:algorithm8}--\eqref{eq:algorithm9} ,
 \eqref{eq:algorithm10}--\eqref{eq:algorithm11} and \eqref{eq:algorithm16}--\eqref{eq:algorithm17} can be easily solved by a preconditioned conjugate gradient algorithm which is detailed discussed in Glowinski\cite{Handbook}. The finite element discretization of \eqref{eq:algorithm4}--\eqref{eq:algorithm5}, \eqref{eq:algorithm6}--\eqref{eq:algorithm7},
\eqref{eq:algorithm12}--\eqref{eq:algorithm13} and \eqref{eq:algorithm14}--\eqref{eq:algorithm15} also have been  discribed in \cite{He-slip-Navier}.

\section{Numerical experiments}
\subsection{Jeffery orbit}
 In Jeffery \cite{Jeffery}, an elliptic particle was positioned in a symmetric domain filled with simple shear rate $\dot{\gamma}$. For Stokes flows, i.e. in the absence of inertia, assuming there is no-slip between particle and fluid, the dynamic equation of Jeffery orbit is given by
\begin{align}
&  \omega = \frac{\dot{\gamma}}{2}\left(-1+p\cos2\Theta\right), \label{eq:Jeffery}
\end{align}
therefore the elliptic particle undergoes a periodic tumbling, where $\omega$ is the angular velocity of the ellipse, $\Theta$ is the angular  orientation, $p$ is a parameter measuring the anisotropy of the ellipse, i.e. $p =\frac{1-e^{2}}{1+e^{2}}$, where $e$ is the aspect ratio defined by $e = \frac{b}{a} \leq 1$, $a$  and $b$  are the major radius and minor radius of the elliptic particle respectively. Since the particle is positioned at the center of a symmetric domain, it has zero translational velocity. 

 We now use the method described in preceding sections to numerically simulate the interactions  between fluid and rigid body. 
  The ellipse rigid particle, which  is located at the domain center $(1,1)$, with $a = \frac{1}{4}$ and $b = \frac{1}{8}$. The computational domain is $[0,2]\times[0,2]$. The solid walls move with speed $W = 2$, in $+x$ direction on the top and $-x$ direction on the bottom, and the particle rotates with an angular velocity $\omega$ which is solved by the numerical scheme introduced in the above sections. The theoretic shear rate $\dot{\gamma} = 2$ and the anisotropy parameter $p = 0.6$ can be calculated from the given parameters.  The external force is $\mathbf{g} = 0$ and we choose $\rho_{f} = 1$, the density of the ellipse particle is the same as the fluid  and $\mu = 1.0$, instead of Navier--Stokes equations, we numerically simulate Stokes equations. The first issue is to verify that $\mathbf{C}_{1}$ and $C_{2}$ will approach to $0$  when $\Delta t$ and $h_{1}$ converge to $0$, which is due to that these two terms are correction terms for translation velocity $\mathbf{U}$ and angular velocity $\omega$. The values of $\mathbf{C}_{1}$ and $C_{2}$ for different $\Delta t$ are shown in Table \ref{Table:test1}, where $\Delta t = O(h^{2}_{1})$. As suggested by \cite{He1} and \cite{He-slip-Navier}, we take $h_{2} = h_{1}$ in our numerical simulations. 
 
\begin{figure}[hp]
\centering
\includegraphics[width=2.5in]{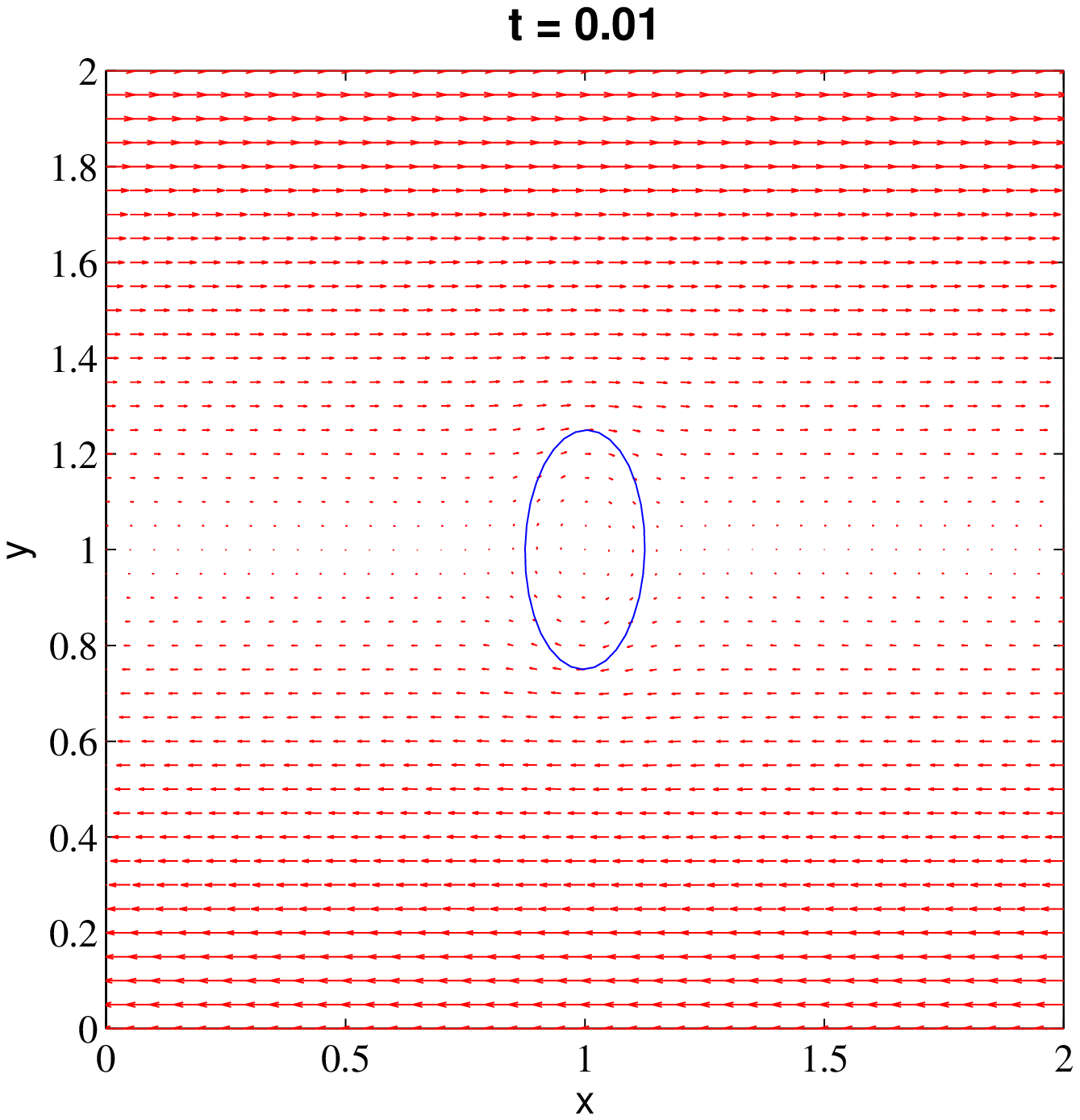}
\includegraphics[width=2.5in]{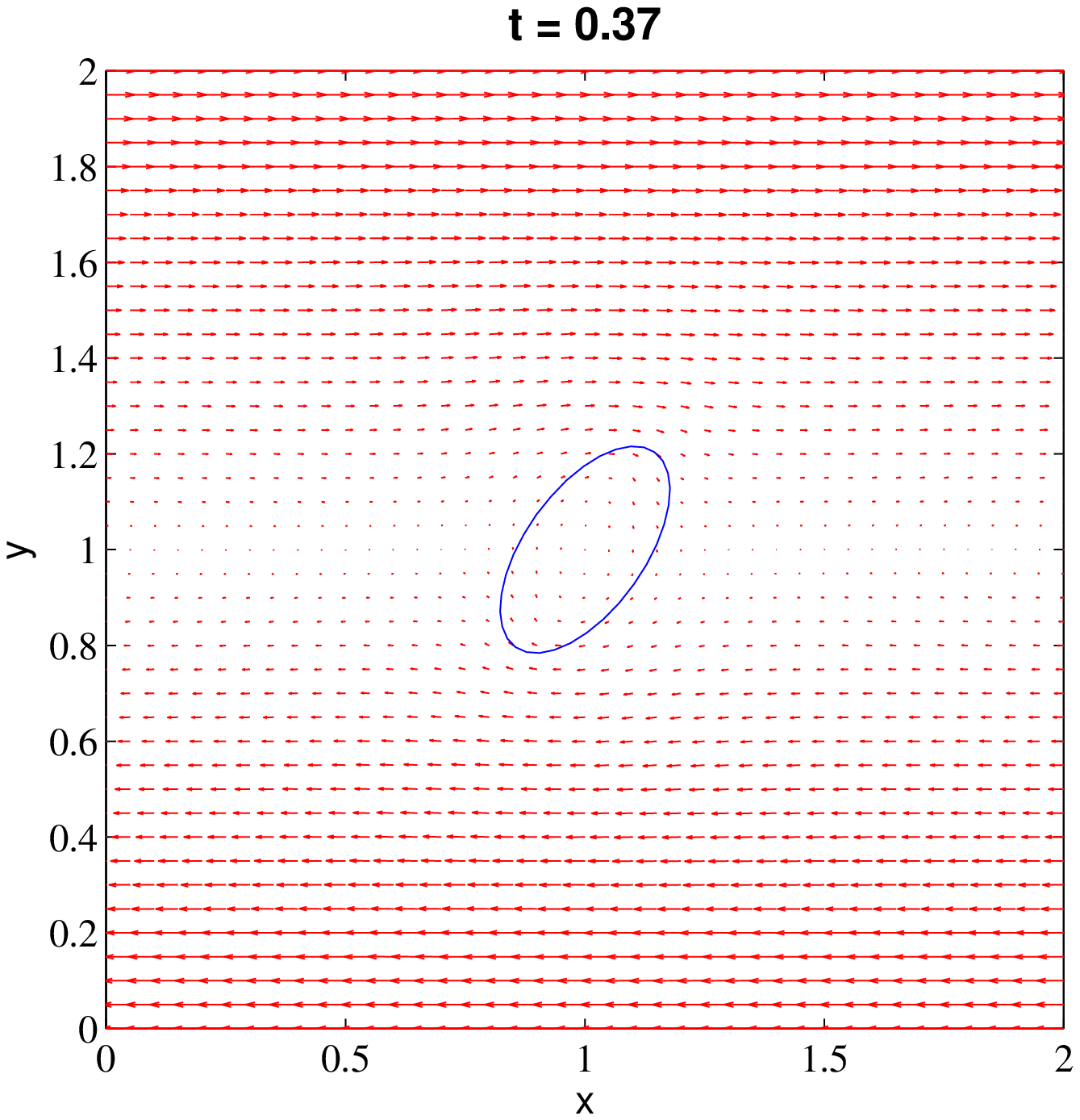}\\
\includegraphics[width=2.5in]{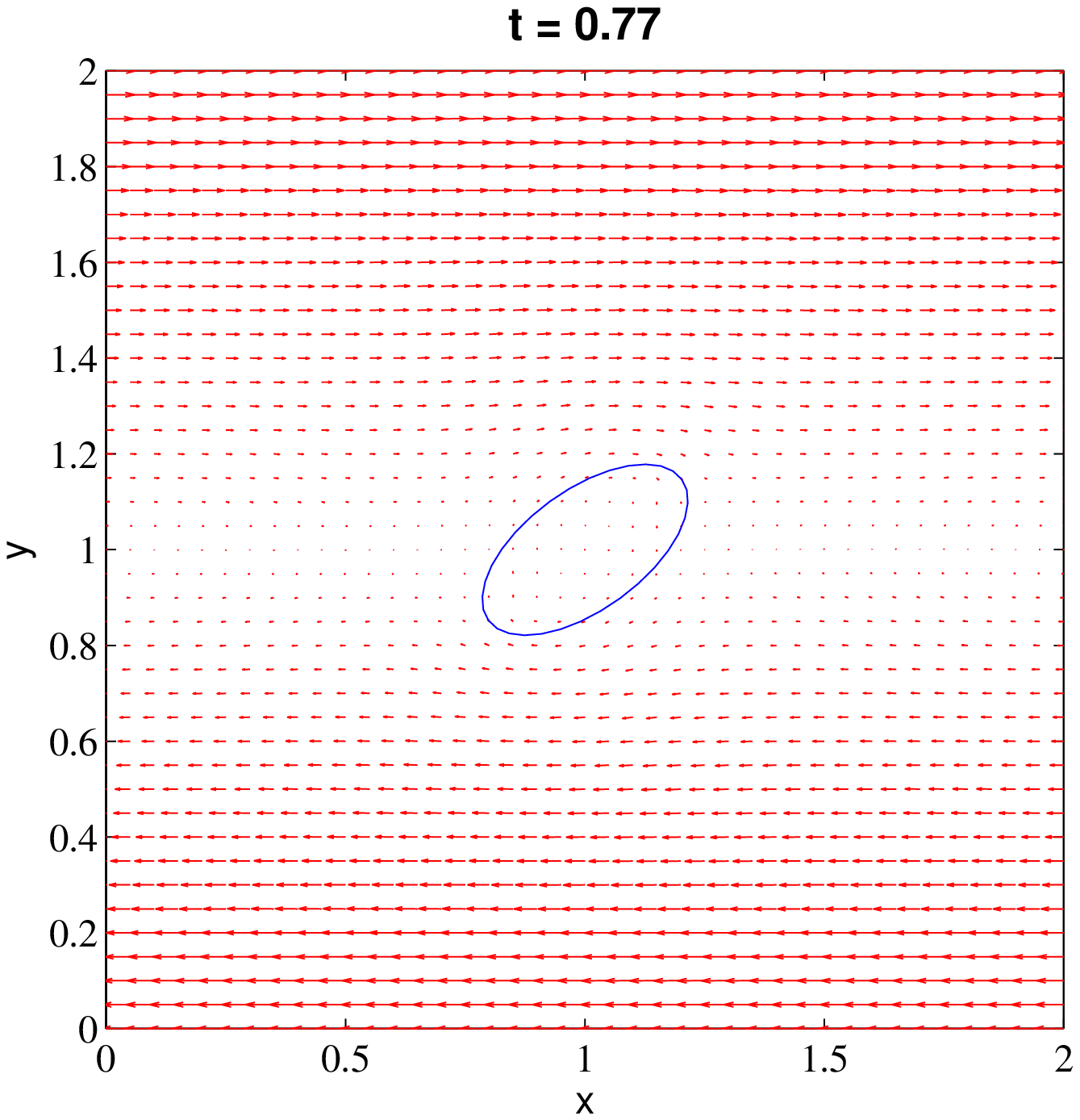}
\includegraphics[width=2.5in]{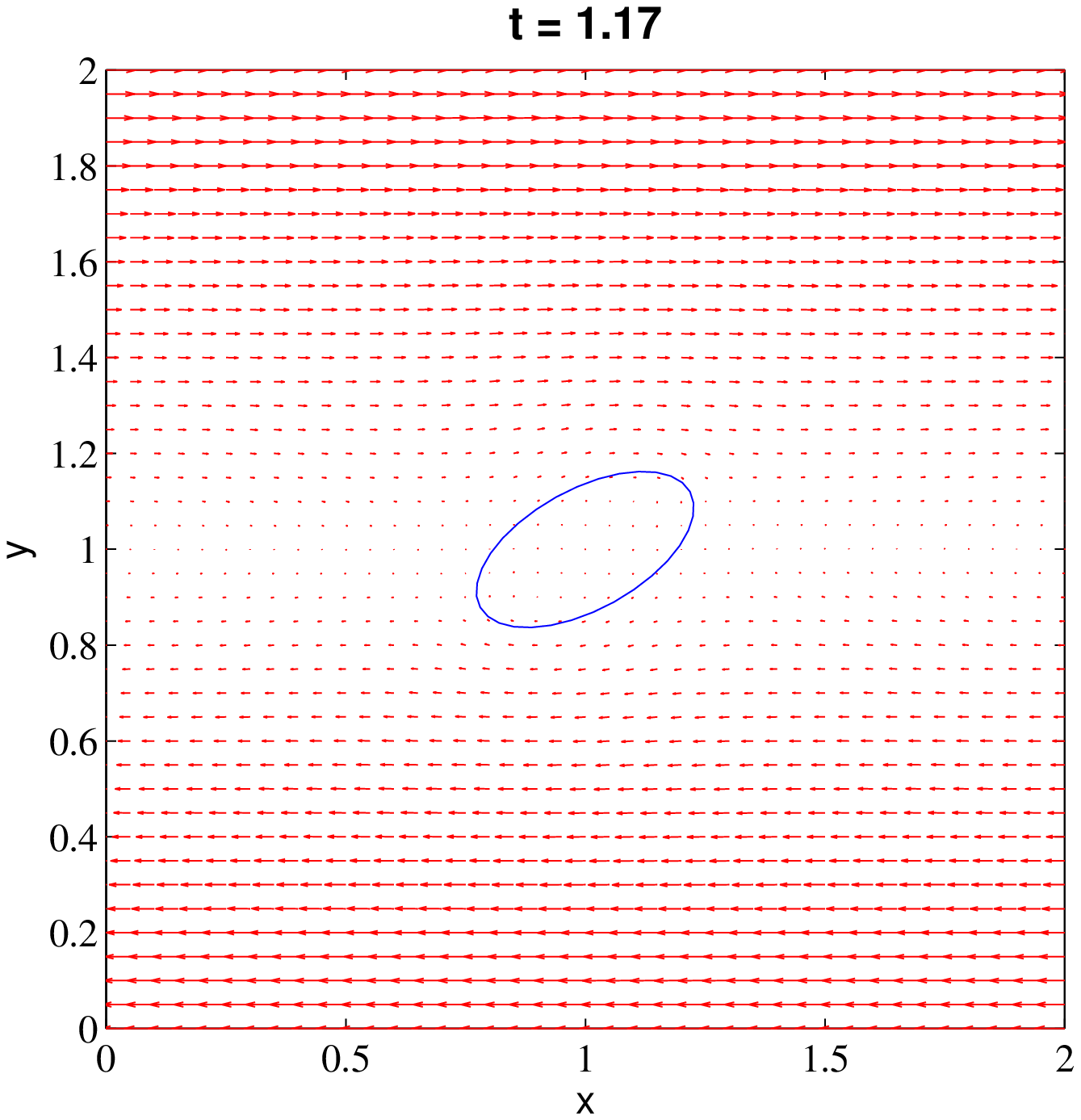}
\vskip -0.2cm
\caption{Snapshots of velocity field  when $t = 0.01, \ 0.37, \  0.77$ and $1.17$ for  $l_{s} = 1/20$.} \label{fig:snapshot}
\end{figure}

\begin{figure}[p]
\centering
\includegraphics[width=9.5cm]{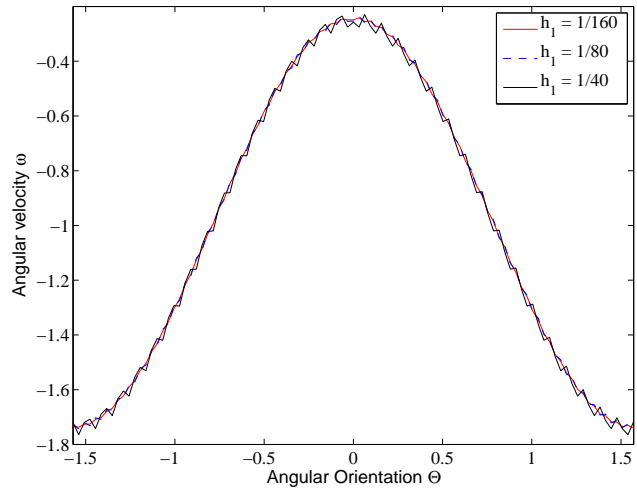}
\vskip -0.2cm
\caption{Comparison of numerical results of different mesh size $h_{1}$ for a fixed $h_{2}$ with $l_{s} = 1/20$.} \label{fig:comparisonmesh}
\end{figure}

\begin{figure}[p]
\centering
\includegraphics[width=9.5cm]{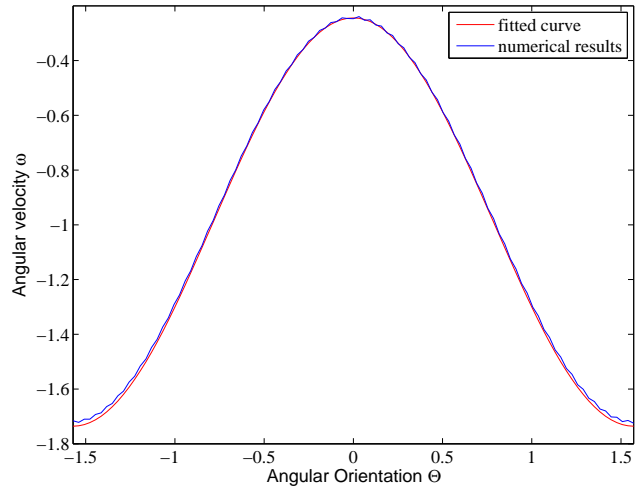}
\vskip -0.2cm
\caption{The angular velocity $\omega$ is plotted as a function of the ellipse orientation $\Theta$ for  $l_{s} = 1/20$.} \label{fig:comparisonfitted}
\end{figure}

\begin{figure}[p]
\centering
\includegraphics[width=9.5cm]{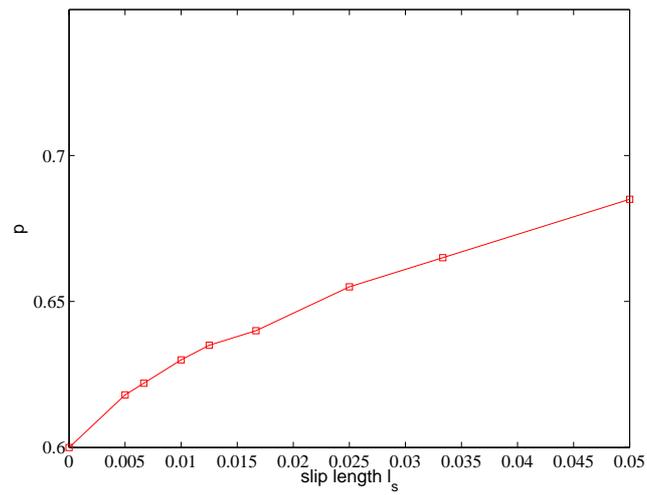}
\vskip -0.2cm
\caption{Numerical solutions of $p$ versus to slip length $l_{s}$.} \label{fig:comparisonp-ls}
\end{figure}

\begin{figure}[p]
\centering
\includegraphics[width=9.5cm]{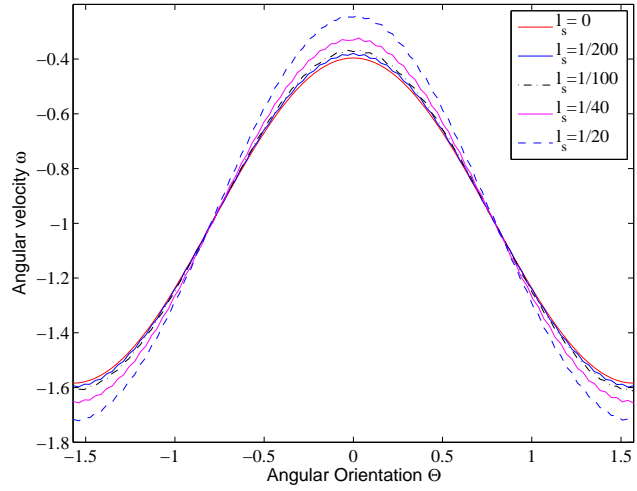}
\vskip -0.2cm
\caption{Numerical solutions of $\omega$ versus to $\Theta$ with $l_{s}= 0, \ 1/200, \  1/100,\  1/40, \ 1/20$.} \label{fig:comparisonls}
\end{figure}

  We will address how  angular velocity $\omega$ changes  with the slip length $l_{s}$. We will numerically verify the effective anisotropy of the particle could be enhanced by the boundary slip, which is presented in \cite{Qian}. Figure \ref{fig:snapshot} shows snapshots of velocity field in different time, and the rotation is observed and translation is almost zero. Next,  we will verify the Jeffery orbit formulation \eqref{eq:Jeffery}.  Using our numerical schemes with slip length $l_{s} = \frac{1}{20}$, we see the dependence of the angular velocity  $\omega$ and the angular  orientation $\Theta$ in Figure \ref{fig:comparisonmesh} with different mesh size, which shows the convergence of our numerical method. The numerical errors (oscillations) can be reduced by mesh refinement. In order to further investigate the Jeffery obit theory with boundary slip, we plot the fitting curve according to \eqref{eq:Jeffery}  and our numerical results in Figure \ref{fig:comparisonfitted}, where the fitted shear rate $ \dot{\gamma} = 1.97$ (close to theoretical value $\dot{\gamma} = 2$), the fitted anisotropy parameter $p = 0.69$ , it is observed that our numerical results show a good cosine curve. In Figure \ref{fig:comparisonp-ls}, we see that as slip length increases, the anisotropy parameter $p$ is also increased, which has been explained in \cite{Qian} that the effective anisotropy of the particle could be enhanced by the boundary slip. The numerical solution of the anisotropy parameter $p$ is computed by equation \eqref{eq:Jeffery}, since numerical solution of $\omega$ can be obtained from our numerical scheme.  With different slip length $l_{s} = 0,\  1/200, \  1/100,\  1/40,  \ 1/20$, we plot the numerical results of $\omega$ along with angular orientation $\Theta$ in Figure \ref{fig:comparisonls}. It is clear to see that our numerical method can resolve the effect of slip length  to the solution efficiently. It is proposed that the effective
  shape of the particle has the major and minor axes given
  by $a -l_{s}$ and $b - l_{s}$ in \cite{Qian}. It is  observed that the slip length affect the amplitude of trigonometric function, which can be obtained by varying the value of axial ratio $p$. The results are shown in Figure \ref{fig:comparisonls}, which  are fairly matched with the that in \cite{Qian}.
  
  \begin{table}[bht]
  	\begin{center}
  		\caption{ \label{Table:test1} The variations of $ \mathbf{C}_{1}$ and $C_{2}$.}
  		\begin{tabular}{ c c c c  }
  			\hline & $\Delta t$  & $\mathbf{C}_{1}$  &  $C_{2}$    \\
  			\hline
  			  & 2.0e-3 &  1.85e-3,\  2.12e-3  & 2.68e-3     \\
  			& 4.88e-4 & 4.26e-4,\ 6.35e-4  & 8.76e-4 \\
  		    & 1.22e-4&  7.32e-5,\ 8.11e-5 &  1.05e-4  \\
  		    & 3.51e-5&  8.34e-6,\ 9.23e-6 &  3.24e-5  \\
  			\hline
  		\end{tabular}
  	\end{center}
  \end{table}



\subsection{Sedimentation of an elliptic particle in a viscous flow}
\begin{figure}[htbp]
	\centering  
	\includegraphics[width=9.5cm]{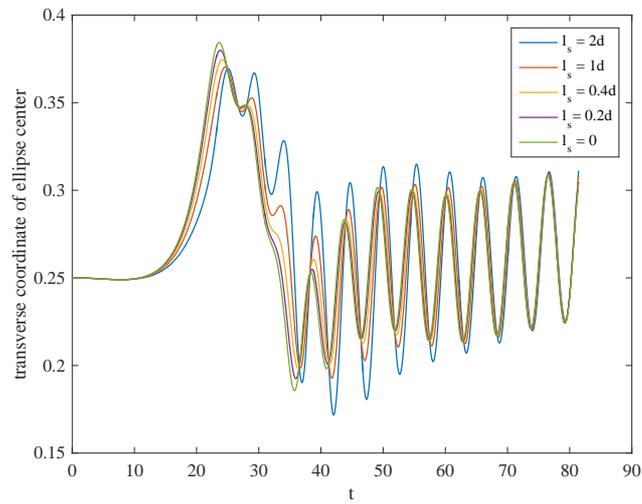}   
	\caption{Transverse coordinate of the ellipse center for $l_s = 2d, \ 1d, \ 0.4d,\ 0.2d, 0$, \ $d = 2a$. } \label{fig:6}
\end{figure}

\begin{figure}[htbp]
	\centering  
	\includegraphics[width=9.5cm]{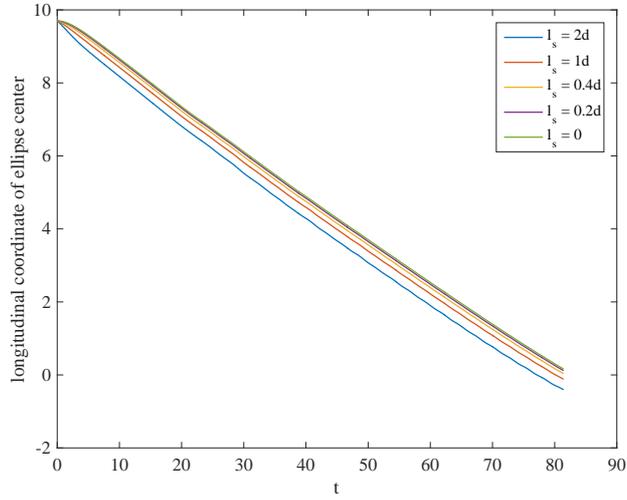}   
	\caption{Longitudinal coordinate of the ellipse center for $l_s = 2d, \ 1d, \ 0.4d,\ 0.2d, 0$, \ $d = 2a$.} \label{fig:7}
\end{figure}

\begin{figure}[htbp]
	\centering  
	\includegraphics[width=9.5cm]{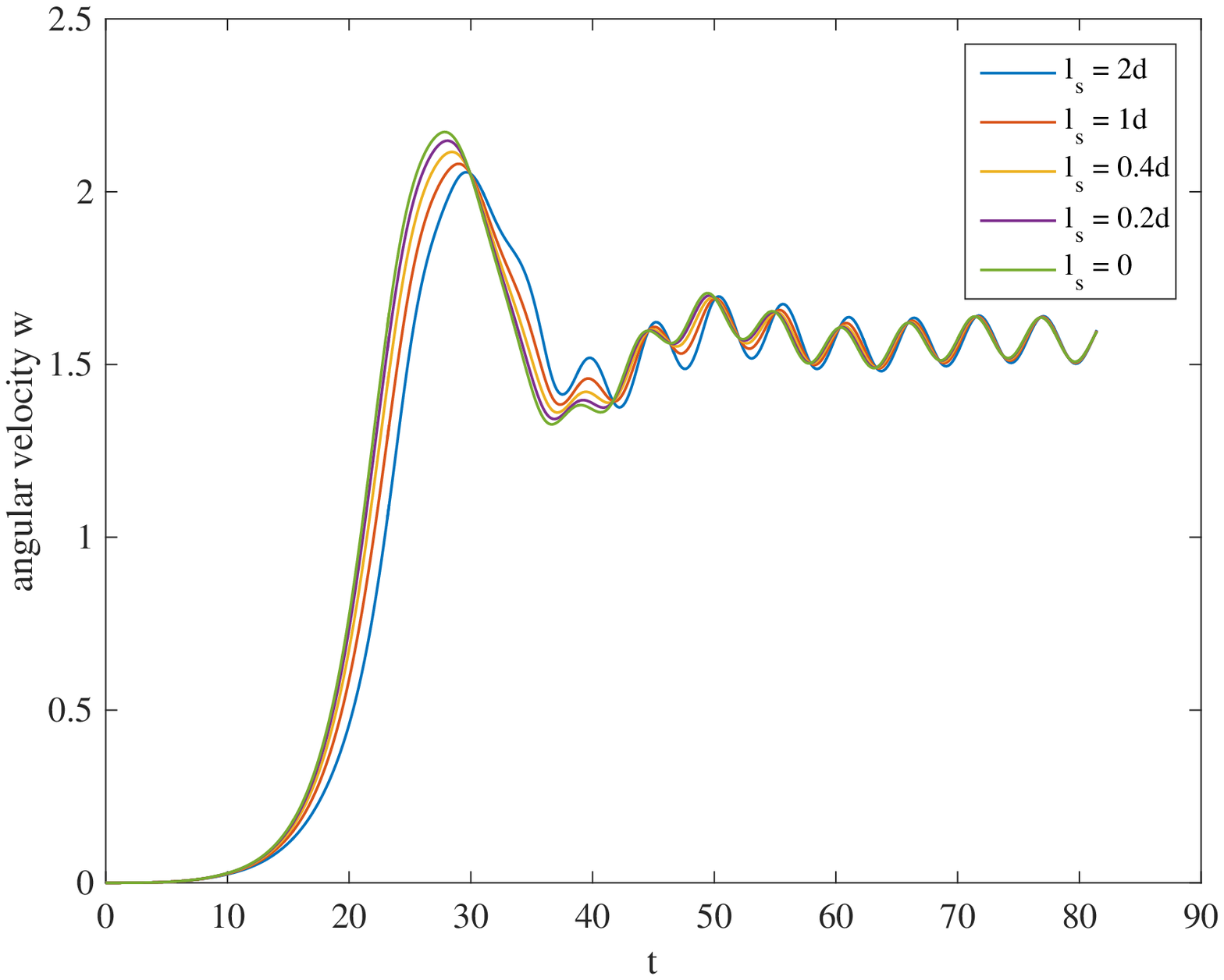}   
	\caption{Angular velocity the ellipse center for $l_s = 2d, \ 1d, \ 0.4d,\ 0.2d, 0$, \ $d = 2a$.} \label{fig:8}
\end{figure}

We consider  the sedimentation of an elliptic particle in a rectangular domain. The ellipse particle is initially static and falls down driven by the gravity $\mathbf{g} = (0,-9.8)^{T}$ and $\mu = 0.1$. The $\rho_f = 1.0$ and $\rho_s = 1.01$ are the density of fluid and solid, respectively. The computational domain $\Omega=[0, 0.5]\times[0,10]$. The elliptic particle is initially located at coordinate $(0.25,9.7)$. The major radius $a$ and minor radius $b$ of the elliptic particle are $0.0625$ and $0.03125$, respectively. From the physical experiments and numerical simulations, we know that the 
particle will undergo a motion called ``drafting and tumbling", which was first numerically demonstrated in \cite{FHJ1994}. The fully incompresible Navier--Stokes are considered. We run the simulations with different slip lengths $ls = 2d, d, 0.4d, 0.2d,0$, where $d = 2a$. The instantaneous transverse coordinate and longitudinal coordinate of the center of the ellipse are shown in Figure \ref{fig:6}--\ref{fig:7}. The angular velocity is shown in Figure \ref{fig:8}. From those results, we clearly find that the particle with different boundary conditions initially perform differential dynamic process 
and finally reach similarly pseudo-steady state. It is observed that slip length will increase the drag coefficient, which can be verified by the numerical results obtained by immersed boundary--lattice Boltzmann method in \cite{WHH2020}. Therefore, the particle falls more faster with increased slip length.


\section{Conclusion}
In this article, we  develop a  least--squares/fictitious domain method for directly numerical simulation of fluid particle interaction with Navier slip boundary condition at the
fluid--particle interface. The method relies on a least--squares formulation of the virtual control type, making it  well--suited for solution by a conjugate gradient algorithm operating in a well--chosen control space.  Numerical experiments show  that our method  works well for anisotropic particle in viscous shear flow with Navier slip boundary condition at the particle surface, showing that the boundary slip can effectively enhance  the anisotropy of the particle, which fairly matches with the results of \cite{Qian}. In the above sections we have been assuming that $\Omega$ contains only one particle. However, we can easily consider cases where $\Omega$ contains a large number particles. In the future, we will generalize our method to the simulation of multiphase flow and particles interactions with boundary slip. 
\vskip 1cm

\section*{Acknowledgments}
The authors would like to thank Prof. Roland Glowinski and Prof. Xiao-Ping Wang for their useful discussions. This research is supported part by National Key R \& D Program of China (2018YFC0830300) and  the National Natural Science Foundation of China (No.11971020). 

\vskip 1.0cm

\section*{References}
\bibliography{journal}
\bibliographystyle{model1b-num-names}

\end{document}